\documentstyle[12pt]{article}

\newcommand{\sect}[1]{\setcounter{equation}{0}\section{#1}}

\newcommand{\app}[1]{\setcounter{section}{0}
\setcounter{equation}{0} \renewcommand{\thesection}{\Alph{section}}
\section{#1}}

\newcommand{\eq}{\begin{equation}}
\newcommand{\eqa}{\begin{eqnarray}}
\newcommand{\en}{\end{equation}}
\newcommand{\ena}{\end{eqnarray}}
\newcommand{\enn}{\nonumber \end{equation}}

\def\sk{\vskip .4cm}
\def\noi{\noindent}
\def\om{\omega}
\def\Om{\Omega}
\def\al{\alpha}
\def\la{\lambda}
\def\be{\beta}
\def\ga{\gamma}
\def\Ga{\Gamma}
\def\del{\delta}
\def\linv{{1 \over \lambda}}
\def\rinv{{1\over {r-r^{-1}}}}

\def\Cb{\bar{C}}

\def\fb{\bar{f}}

\def\rhop{{\rho}^{\prime}}

\def\thetap{{\theta}^{\prime}}

\def\onehalf{{1 \over 2}}

\def\epsi{\varepsilon}
\def\we{\wedge}

\def\de{\delta}

\def\part{\partial}

\def\R#1#2{ R^{#1}_{~~~#2} }
\def\PA#1#2{ P^{#1}_{A~~#2} }

\def\PI#1#2{(P_I)^{#1}_{~~~#2} }
\def\PJ#1#2{ (P_J)^{#1}_{~~~#2} }

\def\Rp#1#2{ (R^+)^{#1}_{~~~#2} }

\def\Rm#1#2{ (R^-)^{#1}_{~~~#2} }
\def\Rinv#1#2{ (R^{-1})^{#1}_{~~~#2} }

\def\Rpm#1#2{(R^{\pm})^{#1}_{~~~#2} }

\def\Rb{{\bf \mbox{\boldmath $R$}}}
\def\Rbo{{\bf \mbox{\boldmath $R$}}}

\def\Rh{{\hat R}}

\def\Rhat#1#2{ \Rh^{#1}_{~~~#2} }

\def\L#1#2{ \La^{#1}_{~~~#2} }

\def\Rhatinv#1#2{ (\Rh^{-1})^{#1}_{~~~#2} }

\def\Z#1#2{ Z^{#1}_{~~~#2} }
\def\X#1#2{ X^{#1}_{~~~#2} }

\def\La{\Lambda}

\def\ff#1#2#3{f_{#1~~~#3}^{~#2}}
\def\MM#1#2#3{M^{#1~~~#3}_{~#2}}
\def\MMc#1#2#3{{M_{\!-}}^{#1~~~#3}_{~#2}}
\def\MMcc#1#2{M_{\!-}{}_{#1}{}^{#2}}
\def\cchi#1#2{\chi^{#1}_{~#2}}
\def\chil#1{\chi_{{}_{#1}}}
\def\ome#1#2{\om_{#1}^{~#2}}
\def\Ome#1#2{\Omega_{#1}^{~#2}}
\def\RRhat#1#2#3#4#5#6#7#8{\La^{~#2~#4}_{#1~#3}|^{#5~#7}_{~#6~#8}}
\def\RRhatinv#1#2#3#4#5#6#7#8{(\La^{-1})^
{~#2~#4}_{#1~#3}|^{#5~#7}_{~#6~#8}}
\def\LL#1#2#3#4#5#6#7#8{\La^{~#2~#4}_{#1~#3}|^{#5~#7}_{~#6~#8}}

\def\Cb{\bf \mbox{\boldmath $C$}}
\def\CC#1#2#3#4#5#6{{\Cb}_{~#2~#4}^{#1~#3}|_{#5}^{~#6}}
\def\cc#1#2#3#4#5#6{C_{~#2~#4}^{#1~#3}|_{#5}^{~#6}}
\def\PIJ#1#2#3#4#5#6#7#8{(P_I,P_J)^{~#2~#4}_{#1~#3}|^{#5~#7}_{~#6~#8}}

\def\ZZ#1#2#3#4#5#6#7#8{Z^{~#2~#4}_{#1~#3}|^{#5~#7}_{~#6~#8}}

\def\C#1#2{ {\bf \mbox{\boldmath $C$}}_{#1}^{~~~#2} }
\def\c#1#2{ C_{#1}^{~~~#2} }

\def\DR{\Delta_R}
\def\DL{\Delta_L}
\def\f#1#2{ f^{#1}_{~~#2} }

\def\T#1#2{ T^{#1}_{~~#2} }

\def\M#1#2{ M_{#1}^{~#2} }
\def\Mc#1#2{ {M_{\!-}}_{#1}^{~#2} }

\def\rminus{r^{-1}}

\def\D{\Delta}

\def\Dp{\Delta^{\prime}}

\def\ep{\epsi^{\prime}}
\def\kp{\kappa^{\prime}}
\def\kpm{\kappa^{\prime -1}}
\def\kpsq{\kappa^{\prime 2}}
\def\km{\kappa^{-1}}

\def\rone{r \rightarrow 1}

\def\IGLqrN{IGL_{q,r}(N)}

\def\SOqrNt{SO_{q,r}(N+2)}
\def\SpqrNt{Sp_{q,r}(N+2)}

\def\ISOqrN{ISO_{q,r}(N)}
\def\ISpqrN{ISp_{q,r}(N)}
\def\ISOqroN{ISO_{q,r=1}(N)}

\def\SqrN{S_{q,r}(N)}
\def\SqrNtwo{S_{q,r}(N+2)}
\def\USqrNtwo{U(S_{q,r}(N+2))}
\def\UISqrN{U(IS_{q,r}(N))}
\def\ISqrN{IS_{q,r}(N)}
\def\SqroNt{S_{q,r=1}(N+2)}
\def\SOqroNt{SO_{q,r=1}(N+2)}
\def\SpqroNt{Sp_{q,r=1}(N+2)}
\def\ISqroN{IS_{q,r=1}(N)}
\def\ISOqroN{ISO_{q,r=1}(N)}

\def\SqrNt{S_{q,r}(N+2)}

\def\Tc{{\cal T}}

\def\Dtwo{\Delta_{N+2}}
\def\epsitwo{\epsi_{N+2}}
\def\kappatwo{\kappa_{N+2}}

\def\Lpm#1#2{L^{\pm #1}_{~~~#2}}
\def\Lmp#1#2{L^{\mp#1}_{~~~#2}}
\def\LLpm{L^{\pm}}

\def\LLp{L^{+}}
\def\LLm{L^{-}}
\def\Lp#1#2{L^{+ #1}_{~~~#2}}
\def\Lm#1#2{L^{- #1}_{~~~#2}}
\def\Dcal#1#2{{\cal D}^{#1}_{~#2}}

\def\limrone{\lim_{r \rightarrow 1}}

\def\chit{{\partial}}

\def\n2{{{N+1} \over 2}}
\def\ap{a^{\prime}}
\def\bp{b^{\prime}}
\def\cp{c^{\prime}}
\def\dpr{d^{\prime}}
\def\Dc{{\cal D}}

\def\Ntwo{{N\over 2}}

\def\bu{\bullet}
\def\ci{\circ}

\def\square{{\,\lower0.9pt\vbox{\hrule \hbox{\vrule height 0.2 cm
\hskip 0.2 cm \vrule height 0.2 cm}\hrule}\,}}
\def\sma#1{\mbox{\footnotesize #1}}
\def\Q.E.D.{\rightline{$\Box$}}
\def\vt{\vartheta}

\begin{document}

\begin{titlepage}
\rightline{DFTT-53/95}
\rightline{IFUP-TH 64/95}
\rightline{December 1995}
\vskip 2em
\begin{center}{\bf BICOVARIANT CALCULUS ON TWISTED $ISO(N)$,
QUANTUM POINCARE GROUP
AND QUANTUM MINKOWSKI SPACE}
\\[3em]
Paolo Aschieri \\[1em]
{\sl Scuola Normale Superiore\\
Piazza dei Cavalieri 7,  56100 Pisa \\and\\Istituto Nazionale
di Fisica Nucleare,\\ Sezione di Pisa, Italy}\\[3em]
Leonardo Castellani\\[1em]
{\sl II Facolt\`a di Scienze M.F.N. di Torino, sede di
Alessandria}\\[.5 em]
{\sl Dipartimento di Fisica Teorica\\
and\\Istituto Nazionale di
Fisica Nucleare\\
Via P. Giuria 1, 10125 Torino, Italy.}  \\[3em]
\end{center}

\begin{abstract}

A bicovariant calculus on the twisted inhomogeneous
multiparametric $q$-groups of the $B_n,C_n,D_n$
type, and on the corresponding quantum planes,
is found by means of a projection
from the bicovariant calculus on
$B_{n+1}$, $C_{n+1}$, $D_{n+1}$.
In particular we obtain the bicovariant calculus
on a dilatation-free $q$-Poincar\'e group
$ISO_{q}(3,1)$, and on the corresponding quantum
Minkowski space.

The classical limit of the $B_n,C_n,D_n$ bicovariant calculus is  
discussed in
detail.
\end{abstract}

\vskip 1.5cm

\noi DFTT-53/95

\noi IFUP-TH 64/95~~~~~~~~~~~~~~~~~~~~~~~~~~~~~~~~~~~
{}~~~~~~~~~~~~~~~~~~~~~~~~~~~~~~~~~q-alg/9601006
\vskip .2cm
\noi \hrule
\vskip.2cm
\noi {\small e-mail: aschieri@ux2sns.sns.it,
castellani@to.infn.it}

\end{titlepage}
\newpage
\setcounter{page}{1}

\sect{Introduction}

We present a bicovariant differential calculus on the
inhomogeneous multiparametric quantum groups of
the $B_n,C_n,D_n$ type, and on the corresponding quantum planes.
Our main
motivation being  an exhaustive study of  the differential calculus
on
the quantum Poincar\'e group  found in ref.s \cite{qpoincare1},
we mainly focus our discussion on the orthogonal inhomogeneous
$q$-groups
$ISO_q(N)$.
All the quantities relevant to their differential calculus are
explicitly constructed.  The results are then directly
applied to the $q$-Poincar\'e group $ISO_q(3,1)$, and to the quantum
Minkowski
space that emerges as the quantum coset
$Fun_{q}(ISO (3,1)/SO(3,1))$.
\sk
The technique used in deriving the differential calculus on
$ISO_q(N)$ or $ISp_q(N)$ is based on a projection from the
bicovariant calculus on $ISO_q(N+2)$ or $ISp_q(N+2)$. This
technique was first proposed in \cite{qigl1} and applied to find
the quantum inhomogeneous groups $IGL_q(N)$ and the corresponding
bicovariant calculus. Their
multiparametric extensions were treated in  \cite{qigl2}.
Other references on inhomogeneous $q$-groups
can be found in \cite{inhom}.
\sk
The projection method was then applied to the multiparametric
$SO_{q,r=1}(N+2)$
 to obtain the bicovariant calculus on $ISO_{q,r=1}(N)$, where $r=1$
corresponds to
the ``minimal deformations" , or twistings, with diagonal $R$-matrix
\cite{qpoincarediff1}.
The gauging of the resulting
deformed $q$-Poincar\'e  ``Lie algebra" leads to the $q$-gravity
theory discussed in the same references. The structure of the
multiparametric
inhomogeneous $q$-groups $ISO_{q,r}(N)$ obtained via the
projection technique was studied in detail in ref. \cite{qpoincare1},
where
a dilatation-free $q$-Poincar\'e group depending on one real
deformation
parameter was found. Absence of dilatations requires $r=1$.
\sk
In the present paper the bicovariant calculus on $ISO_{q,r}(N)$ and
$ISp_{q,r}(N)$ with $r=1$ is obtained after a detailed study of the  
homogeneous
orthogonal and symplectic
$q$-groups in the $\rone$ limit. The necessity of taking $r=1$ is  
discussed.
The functionals of the universal enveloping algebra  
$U(SO_{q,r=1}(N+2))
[U(Sp_{q,r=1}(N+2))]$
relevant for the construction of a bicovariant calculus are analized  
and
``projected''
to well defined functionals on $ISO_{q,r=1}(N) [ISp_{q,r=1}(N)]$. The
differential calculus is
found and explicitly formulated in terms of these ``projected''  
functionals,
that in the commutative
limit become the tangent vectors to the inhomogeneous orthogonal
 [symplectic] groups.
In this general setting we are able to retrieve
all the results of \cite{qpoincarediff1}   (where these functionals
were
only given in terms of their action on the adjoint $q$-group
elements) in a
direct way, and to
clarify important points. For example we will easily see how
the typical ``abundance" of left-invariant one-forms can be lifted in
the $r=1$
case.
\sk
In our framework the bicovariant calculus on the orthogonal
multiparametric quantum plane follows almost as a
corollary.
\sk
In Section 2 we briefly review the basics of $B_n,C_n,D_n$
multiparametric
quantum groups, mainly to establish notations.  In Section 3 we
recall the
$R$-matrix formulation of  $ISO_{q,r}(N)$ and $ISp_{q,r}(N)$
of ref. \cite{qpoincare1}, and discuss the real forms
yielding  $ISO_{q,r}(n,n;\Rb)$, $ISO_{q,r}(n,n+1;\Rb)$,
$ISp_{q,r}(n,\Rb)$ and  $ISO_{q,r}(n+1,n-1;\Rb)$, this last
being the one used in \cite{qpoincarediff1,qpoincare1}
to obtain the quantum Poincar\'e group $ISO_{q,r}(3,1)$.
The universal enveloping algebra and the
bicovariant calculus on multiparametric $B_n,C_n,D_n$ $q$-groups
(and their real forms)
are given in Section 4 and 5 respectively. In Section 6 we examine
the case $r=1$.  We clarify
some issues related to the classical limit and
see how in this limit some tangent vectors become
linearly independent, thus providing the correct classical dimension
of the tangent space. A similar mechanism occurs for the  
left-invariant
one-forms.
In  Section 7 the bicovariant calculus on $ISO_{q,r=1}(N)$ and its  
real forms
is constructed. Finally Section 8 deals with the differential
calculus on the
orthogonal quantum plane $Fun_{q,r=1}(ISO(N)/SO(N))$.

\sect{$B_n,C_n,D_n$ multiparametric quantum groups}

The $B_n, C_n, D_n$ multiparametric
quantum groups are freely generated by the noncommuting
matrix elements $\T{a}{b}$ (fundamental representation) and
the identity $I$. The noncommutativity is controlled by the $R$
matrix:
\eq
\R{ab}{ef} \T{e}{c} \T{f}{d} = \T{b}{f} \T{a}{e} \R{ef}{cd}
\label{RTT}
\en
which satisfies the quantum Yang-Baxter equation
\eq
\R{a_1b_1}{a_2b_2} \R{a_2c_1}{a_3c_2} \R{b_2c_2}{b_3c_3}=
\R{b_1c_1}{b_2c_2} \R{a_1c_2}{a_2c_3} \R{a_2b_2}{a_3b_3}, \label{QYB}
\en
a sufficient condition for the consistency of the
``$RTT$" relations (\ref{RTT}).  The $R$-matrix components
$\R{ab}{cd}$
depend
continuously on a (in general complex)
set of  parameters $q_{ab},r$.  For $q_{ab}=r$ we
recover the uniparametric $q$-groups of ref. \cite{FRT}. Then
$q_{ab} \rightarrow 1, r \rightarrow 1$ is the classical limit for
which
$\R{ab}{cd} \rightarrow \de^a_c \de^b_d$ : the
matrix entries $\T{a}{b}$ commute and
become the usual entries of the fundamental representation. The
multiparametric $R$ matrices for the $A,B,C,D$ series can be
found in \cite{Schirrmacher}  (other ref.s on multiparametric
$q$-groups are given in \cite{multiparam1,multiparam2}). For
the $B,C,D$ case they read:
\eq
\begin{array}{ll}
\R{ab}{cd}=&\delta^a_c \delta^b_d [{r\over q_{ab}} +
(r-1) \delta^{ab}+(\rminus - 1)
\de^{a\bp}] (1-\de^{a n_2})
+\de^a_{n_2} \de^b_{n_2} \de^{n_2}_c \de^{n_2}_d  \\
&+(r-r^{-1})
[\theta^{ab} \delta^b_c \de^a_d - \epsilon_a\epsilon_c
\theta^{ac} r^{\rho_a - \rho_c}
\de^{\ap b} \de_{\cp d}]
\end{array}
\label{Rmp}
\en
\noi where $\theta^{ab}=1$ for $a> b$
and $\theta^{ab}=0$ for $ a \le b$; we define
 $n_2 \equiv \n2$ and primed indices as $\ap \equiv N+1-a$. The
indices
run on $N$ values ($N$=dimension of
the fundamental representation $\T{a}{b}$),
with $N=2n+1$ for $B_n  [SO(2n+1)]$,  $N=2n$
for $C_n [Sp(2n)]$, $D_n [SO(2n)]$.
The terms with the index $n_2$ are present
only for the $B_n$ series. The $\epsilon_a$ and
$\rho_a$ vectors are given by:
\eq
\epsilon_a=
\left\{ \begin{array}{ll} +1 & \mbox{for $B_n$, $D_n$} ,\\
                            +1 & \mbox{for $C_n$ and $a \le n$},\\
                              -1  & \mbox{for $C_n$ and $a > n$}.
           \end{array}
                                    \right.
\en
\eq
(\rho_1,...\rho_N)=\left\{ \begin{array}{ll}
         (\Ntwo -1, \Ntwo -2,...,
{1\over 2},0,-{1\over 2},...,-\Ntwo+1)
                   & \mbox{for $B_n$} \\
           (\Ntwo,\Ntwo -1,...1,-1,...,-\Ntwo) & \mbox{for $C_n$} \\
           (\Ntwo -1,\Ntwo -2,...,1,0,0,-1,...,-\Ntwo+1) & \mbox{for
$D_n$}
                                             \end{array}
                                    \right.
\en
Moreover the following relations reduce the number of independent
$q_{ab}$ parameters \cite{Schirrmacher}:
\eq
q_{aa}=r,~~q_{ba}={r^2 \over q_{ab}}; \label{qab1}
\en
\eq
q_{ab}={r^2 \over q_{a\bp}}={r^2 \over q_{\ap b}}=q_{\ap\bp}
 \label{qab2}
\en
\noi where (\ref{qab2}) also implies $q_{a\ap}=r$. Therefore
 the $q_{ab}$ with $a < b \le {N\over 2}$ give all the $q$'s.
\sk
It is useful to list the nonzero complex components
of the $R$ matrix (no sum on repeated indices):
\eqa
& &\R{aa}{aa}=r , ~~~~~~~~~~~~~~\mbox{\footnotesize
 $a \not= {n_2}$ } \cr
& &\R{a\ap}{a\ap}=r^{-1} ,  ~~~~~~~~~~\mbox{\footnotesize
 $a \not= {n_2}$ } \cr
& &\R{{n_2}{n_2}}{{n_2}{n_2}}= 1\cr
& &\R{ab}{ab}={r \over q_{ab}} ,~~~~~~~~~~~~\mbox{\footnotesize
 $a \not= b$,  $\ap \not= b$}\cr
& &\R{ab}{ba}=r-r^{-1} , ~~~~~~~\mbox{\footnotesize
$a>b, \ap \not= b $}\cr
& &\R{a\ap}{\ap a}=(r-r^{-1})(1-\epsilon r^{\rho_a-\rho_{\ap}}) ,
{}~~~~~~\mbox{\footnotesize
$a>\ap $}\cr
& &\R{a\ap}{b \bp}=-(r-r^{-1})\epsilon_a\epsilon_b
r^{\rho_a-\rho_b} , ~~~~~~~
\mbox{\footnotesize $~~a>b ,~ \ap \not= b $}
\label{Rnonzero}
\ena
where $\epsilon=\epsilon_a \epsilon_{\ap}$, i.e.
 $\epsilon=1$ for $B_n$, $D_n$ and $\epsilon=-1$
for $C_n$.
\sk
{\sl Remark 2.1 :}  The matrix $R$ is upper triangular, that is
$\R{ab}{cd}=0$ if  [$\sma{$a=c$}$ and $\sma{$b<d$}$] or
$\sma{$a<c$}$,
and has the following properties:
\eq
R^{-1}_{q,r}=R_{q^{-1},r^{-1}}~~;~~~
(R_{\!q,r})^{ab}{}_{cd}=(R_{\!q,r})^{\cp\dpr}{}_{\ap\bp}~~;~~~
(R_{q,r})^{ab}{}_{cd}=(R_{\!p,r})^{dc}{}_{ba}
\label{Rprop1}
\en
where
$q,r$ denote  the set of parameters $q_{ab},r$, and $p_{ab}\equiv
q_{ba}$.
\sk
The inverse $R^{-1}$ is defined by
$\Rinv{ab}{cd} \R{cd}{ef}=\de^a_e \de^b_f=\R{ab}{cd}
\Rinv{cd}{ef}$.
Eq. (\ref{Rprop1}) implies
that for $|q|=|r|=1$, ${\bar R}=R^{-1}$.
\sk
{\sl Remark 2.2:} Let $R_r$ be the uniparametric $R$ matrix
for the $B, C, D$ q-groups. The multiparametric $R_{q,r}$
matrix is obtained from $R_r$ via the transformation
\cite{multiparam1,Schirrmacher}
\eq
R_{q,r}=F^{-1}R_rF^{-1}
\en
where $(F^{-1})^{ab}_{~~cd}$ is a diagonal matrix
in the index couples $ab$, $cd$:
\eq
F^{-1}\equiv diag (\sqrt{{r \over q_{11}}} ,
\sqrt{{r \over q_{12}}} , ... ~ \sqrt{{r \over q_{NN}}})
\label{effe}
\en
\noi and $ab$, $cd$ are ordered as
in the $R$ matrix.
Since $\sqrt{{r \over q_{ab}}} =(\sqrt{{ q_{ba}\over r}})^{-1}$
and $q_{a\ap}=q_{b\bp}$, the non diagonal
elements of $R_{q,r}$ coincide with those of $R_r$.
The matrix $F$ satisfies $F_{12}F_{21}=1$ i.e.
$F^{ab}{}_{ef}F^{fe}{}_{dc}=\delta^a_c\delta^b_d $,
the quantum Yang-Baxter equation $F_{12}F_{13}F_{23}
=F_{23}F_{13}F_{12}$ and the relations
$(R_r)_{12}F_{13}F_{23}=F_{23}F_{13}(R_r)_{12}$.
Note that for $r=1$ the multiparametric $R$ matrix reduces to
$R=F^{-2}$.
\sk
{\sl Remark 2.3:} Let $\Rh$ the matrix defined by
$\Rhat{ab}{cd} \equiv \R{ba}{cd}$.  Then the multiparametric
$\Rh_{q,r}$ is
obtained from $\Rh_r$ via the similarity transformation
\eq
\Rh_{q,r}=F\Rh_rF^{-1}
\en
The characteristic equation and the projector decomposition
of $\Rh_{q,r}$ are therefore the same as in the uniparametric case:
\eq
(\Rh-rI)(\Rh+r^{-1}I)(\Rh-\epsilon r^{\epsilon-N} I)=0 \label{cubic}
\en
\eq
\Rh-\Rh^{-1}=(r-r^{-1}) (I-K) \label{extrarelation}
\en
\eq
\Rh=r P_S - r^{-1} P_A+\epsilon r^{\epsilon-N}P_0  \label{RprojBCD}
\en
with
\eq
\begin{array}{ll}
&P_S={1 \over {r+\rminus}} [\Rh+\rminus I-(\rminus+\epsilon
r^{\epsilon-N})P_0]\\
&P_A={1 \over {r+\rminus}} [-\Rh+rI-(r-\epsilon r^{\epsilon-N})P_0]\\
&P_0= Q_N(r) K\\
&Q_N(r) \equiv (C_{ab} C^{ab})^{-1}={{1-r^2} \over
{(1-\epsilon r^{N+1-\epsilon})(1+\epsilon r^{-N+1+\epsilon})}}~,~~~~
K^{ab}_{~~cd}=C^{ab} C_{cd}\\
&I=P_S+P_A+P_0
\end{array}
\label{projBCD}
\en
To prove (\ref{extrarelation}) in the multiparametric case note that
$F_{12} K_{12} F_{12}^{-1}=K_{12}$.
Orthogonality  (and symplecticity) conditions can be
imposed on the elements $\T{a}{b}$, consistently
with  the $RTT$ relations (\ref{RTT}):
\eq
C^{bc} \T{a}{b}  \T{d}{c}= C^{ad},~~~
C_{ac} \T{a}{b}  \T{c}{d}=C_{bd} I \label{Torthogonality}
\en
\noi where the (antidiagonal) metric is :
\eq
C_{ab}=\epsilon_a r^{-\rho_a} \de_{a\bp} \label{metric}
\en
\noi and its inverse $C^{ab}$
satisfies $C^{ab} C_{bc}=\de^a_c=C_{cb} C^{ba}$.
We see
that for the orthogonal series, the matrix elements of the metric
and the inverse metric coincide,
while for the symplectic series there is a change of sign:
$C^{ab}=\epsilon C_{ab}$. Notice also the symmetry
$C_{ab}=C_{\bp\ap}$.

The consistency of (\ref{Torthogonality}) with the $RTT$ relations
is due to the identities:
\eq
C_{ab} \Rhat{bc}{de} = \Rhatinv{cf}{ad} C_{fe} ,~~~
 \Rhat{bc}{de} C^{ea}=C^{bf} \Rhatinv{ca}{fd} \label{crc}
\en
\noi These identities
 hold also for $\Rh \rightarrow \Rh^{-1}$ and can be proved using
the explicit expression (\ref{Rnonzero}) of $R$.
\sk

We note the useful relations, easily deduced from (\ref{RprojBCD}):
\eq
C_{ab}\Rhat{ab}{cd}=\epsilon r^{\epsilon-N}C_{cd} ,~~~
C^{cd}\Rhat{ab}{cd}=\epsilon r^{\epsilon-N}C^{ab}  \label{CR}
\en

The co-structures of the $B_n,C_n,D_n$ multiparametric quantum
groups have the same form as in the uniparametric case:
the coproduct
$\D$, the counit $\epsi$ and the coinverse $\kappa$ are given by
\eqa
& & \D(\T{a}{b})=\T{a}{b} \otimes \T{b}{c}  \label{cos1} \\
& & \epsi (\T{a}{b})=\delta^a_b\\
& & \kappa(\T{a}{b})=C^{ac} \T{d}{c} C_{db}
\label{cos2}
\ena
Four conjugations (i.e. algebra antihomomorphism, coalgebra
homomorphism
and involution, satisfying $\kappa(\kappa(T^*)^*)=T$)
can be defined, but only two of these can be extended
to the corresponding inhomogeneous groups \cite{qpoincare1}.
These two are defined as follows:
\sk
$\bullet$~~trivially as $T^*=T$. Compatibility with the
$RTT$ relations (\ref{RTT}) requires ${\bar R}_{q,r}=R^{-1}_{q,r}=
R_{q^{-1},r^{-1}}$,
i.e. $|q|=|r|=1$. Then the $CTT$ relations are invariant under
$*$-conjugation.  The corresponding real forms are
$SO_{q,r}(n,n;\Rbo)$, $SO_{q,r}(n,n+1;\Rbo)$
(for N even and odd respectively) and $Sp_{q,r}(n;\Rbo)$.
\sk
$\bullet$~~
on the orthogonal
quantum groups $SO_{q,r}(2n,{\bf C})$,  extending to the
multiparametric case the one
proposed by
the authors of ref. \cite{Firenze1} for $SO_{q}(2n,{\bf C})$.
The conjugation is defined by:
\eq
(\T{a}{b})^*={\cal D}^a_{~c} \T{c}{d}
{\cal D}^d_{~b} \label{Tconjugation}
\en
\noi ${\cal D}$ being the matrix that
exchanges the index $n$ with the index $n+1$.
This conjugation is compatible with the coproduct: $\D (T^*)=(\D
T)^*$;
for $|r|=1$ it is also  compatible with the
orthogonality relations (\ref{Torthogonality}) (due to
${\bar C}=C^T$ and also $\Dc C \Dc = C$) and with the antipode:
$\kappa(\kappa(T^*)^*)=T$.  Compatibility with the $RTT$ relations
is easily seen to require
\eq
({\bar R})_{n \leftrightarrow n+1}=R^{-1},~~~\mbox{ i.e. }~~{\cal
D}_1 {\cal
D}_2 R_{12}
{\cal D}_1 {\cal D}_2 = {\overline R_{12}^{-1} } \label{Rprop2}
\en
\noi which implies

i) $|q_{ab}|=|r|=1$
for $a$ and $b$ both different from $n$ or $n+1$;

ii) $q_{ab}/r \in {\bf R}$
when at least one of the indices $a,b$ is equal
to $n$ or $n+1$.
\sk
\noi This conjugation leads to the
real form $SO_{q,r}(n+1,n-1;\Rbo)$, and
is in fact  the one  needed  to obtain $ISO_{q,r}(3,1;\Rbo)$, as
 discussed  in ref.s \cite{qpoincarediff1,qpoincare1} and later in
this paper.

\sk
{\sl Remark 2.4:}
Using formula (\ref{Rmp}) or (\ref{Rnonzero}),
we find that  the $\R{AB}{CD}$  matrix for the
$SO_{q,r}(N+2)$ and $Sp_{q,r}(N+2)$ quantum groups
can be decomposed
in terms of  $SO_{q,r}(N)$ and $Sp_{q,r}(N)$ quantities
as follows (splitting the index {\small A} as
{\small A}=$(\circ, a, \bullet)$, with $a=1,...N$):
\eq
\R{AB}{CD}=\left(  \begin{array}{cccccccccc}
   {}&\circ\circ&\circ\bullet&\bullet
          \circ&\bullet\bullet&\circ d&\bullet d
      &c \circ&c\bullet&cd\\
   \circ\circ&r&0&0&0&0&0&0&0&0\\
   \circ\bullet&0&r^{-1}&0&0&0&0&0&0&0\\
   \bullet\circ&0&f(r)&r^{-1}&0&0&0&0&0&-\epsilon C_{cd} \lambda
r^{-\rho}\\
\bullet\bullet&0&0&0&r&0&0&0&0&0\\
\circ b&0&0&0&0&{r\over q_{\circ b}} \de^b_d&0&0&0&0\\
\bullet b&0&0&0&0&0&{r\over q_{\bullet b} }
\de^b_d&0&\lambda\de^b_c&0\\
a\circ&0&0&0&0&\lambda\de^a_d&0&{r \over q_{a \circ} } \de^a_c&0&0\\
a\bullet&0&0&0&0&0&0&0&{r\over q_{a \bullet}} \de^a_c&0\\
ab&0&-C^{ba} \lambda r^{-\rho}
&0&0&0&0&0&0&\R{ab}{cd}\\
\end{array} \right) \label{Rbig}
\en
\noi where $\R{ab}{cd}$ is the $R$ matrix for  $SO_{q,r}(N)$
or $Sp_{q,r}(N)$, $C_{ab}$ is the corresponding
metric,  $\lambda \equiv r-r^{-1}$,
$\rho={{N+1-\epsilon}\over 2}~(r^{\rho}=C_{\bullet \circ})$
and $f(r) \equiv \lambda (1-\epsilon r^{-2\rho})$.
The sign $\epsilon$ has been defined after eq. s
(\ref{Rnonzero}).
\sk

\sect{The quantum inhomogeneous groups
$\ISOqrN$ and $\ISpqrN$}

An $R$-matrix formulation for the
quantum inhomogeneous groups $\ISOqrN$ and $\ISpqrN$
was obtained in ref. \cite{qpoincare1}, in terms of
the $\R{AB}{CD}$  matrix for the
$SO_{q,r}(N+2)$ and $Sp_{q,r}(N+2)$ quantum groups.
It was found that the quantum inhomogeneous groups
$\ISOqrN$ and $\ISpqrN$ are freely generated by the
non-commuting matrix elements $\T{A}{B}$
 [{\small A}=$(\circ, a, \bullet)$, with $a=1,...N$)] and the
identity
$I$,
modulo the relations:
\eq
\T{a}{\circ}=\T{\bullet}{b}=\T{\bullet}{\circ}=0 , \label{Tprojected}
\en
\noi the $RTT$ relations
\eq
\R{AB}{EF} \T{E}{C} \T{F}{D} = \T{B}{F} \T{A}{E} \R{EF}{CD},
\label{RTTbig}
\en
\noi and the orthogonality (symplecticity) relations
\eq
C^{BC} \T{A}{B}  \T{D}{C}= C^{AD} ,~~~
C_{AC} \T{A}{B}  \T{C}{D}=C_{BD} \label{CTTbig}
\en
The co-structures of $\ISOqrN$ and $\ISpqrN$
are simply given by:
\eq
\D (\T{A}{B})=\T{A}{C} \otimes \T{C}{B},
{}~~\kappa (\T{A}{B})=C^{AC} \T{D}{C} C_{DB} ,~~
\epsi (\T{A}{B})=\de^A_B~. \label{costructuresbig}
\en
After decomposing the indices {\small A}=$(\circ, a, \bullet)$, and
defining:
\eq
u\equiv \T{\circ}{\circ},~~v\equiv
\T{\bullet}{\bullet},~~z\equiv
\T{\circ}{\bullet},~~
x^a \equiv \T{a}{\bullet},~~y_a \equiv \T{\circ}{a} \label{names}
\en
\noi the relations (\ref{RTTbig}) and (\ref{CTTbig}) become
\eqa
& &\R{ab}{ef} \T{e}{c} \T{f}{d} = \T{b}{f} \T{a}{e} \R{ef}{cd}
\label{PRTT11}\\
& &\T{a}{b} C^{bc} \T{d}{c}=C^{ad} I \label{PRTT31}\\
& &\T{a}{b} C_{ac} \T{c}{d} = C_{bd} I \label{PRTT32}
\ena
\eqa
& &\T{b}{d} x^a={r \over q_{d\bullet}} \R{ab}{ef} x^e \T{f}{d}
\label{PRTT33}\\
& &\PA{ab}{cd} x^c x^d=0 \label{PRTT13}\\
& &\T{b}{d} v={q_{b\bullet}\over q_{d\bullet}} v \T{b}{d}\\
& &x^b v=q_{b \bullet} v x^b \label{PRTT15}\\
& & uv=vu=I \label{PRTT21}\\
& &u x^b=q_{b\bullet} x^bu \label{PRTT22}\\
& &u \T{b}{d}={q_{b\bullet}\over q_{d\bullet}} \T{b}{d} u
\label{PRTT24}
\ena
\eq
y_b=-r^{\rho} \T{a}{b} C_{ac} x^c u \label{ipsilon}
\en
\eq
z=-{1\over {(r^{-\rho}+\epsilon r^{\rho-2})}} x^b C_{ba} x^a u
\label{PRTT44}
\en
\noi where $q_{a\bullet}$ are $N$  complex parameters
related by $q_{a\bullet} = r^2 /q_{\ap\bullet}$, with $\ap = N+1-a$.
Note that in the symplectic case, $x^b C_{ba} x^a=0$ so that  the
constraint  (\ref{PRTT44}) reads $z=0$.
The matrix $P_A$ in eq. (\ref{PRTT13}) is the $q$-antisymmetrizer for
the
$B,C,D$ $q$-groups given by (cf.  (\ref{projBCD})):
\eq
\PA{ab}{cd}=- {1 \over {r+\rminus}}
(\Rhat{ab}{cd}-r\de^a_c \de^b_d + {r-r^{-1} \over
\epsilon r^{N-1-\epsilon} +1} C^{ab} C_{cd}). \label{PA}
\en
The last two relations (\ref{ipsilon}) - (\ref{PRTT44})
are constraints, showing that
the $\T{A}{B}$ matrix elements in eq. (\ref{RTTbig})
are really a {\sl redundant} set. This redundance
is necessary if we want to express the $q$-commuations
of the $\ISOqrN$ and $\ISpqrN$ basic group elements
as $RTT=TTR$  (i.e. if we want an $R$-matrix
formulation).  Remark that,
in the $R$-matrix formulation
for $\IGLqrN$, {\sl all}  the $\T{A}{B}$
are independent \cite{qigl1,qigl2}. Here
we can take as independent generators the
elements
\eq
{}~~\T{a}{b} , x^a , v , u\equiv v^{-1}
\mbox{ and the identity }  I~~~~~(a=1,...N)
\en
The co-structures on the $\ISOqrN$ generators can be read from
(\ref{costructuresbig})
after decomposing the indices \sma{$A=\circ,a,\bu$}:
\eqa
& &\D (\T{a}{b})=\T{a}{c} \otimes \T{c}{b}~,~~
\D (x^a)=\T{a}{c} \otimes x^c + x^a \otimes v ~,
\label{Pcoproduct2}\\
& &\D (v)=v \otimes v~,~~\D (u)=u \otimes u~,
\ena
\eqa
& &\kappa (\T{a}{b})
=C^{ac} \T{d}{c} C_{db}=\epsilon_a\epsilon_b
r^{-\rho_a+\rho_{b}}~ \T{\bp}{\ap}~,\\
& &\kappa (x^a)
=- \kappa (\T{a}{c}) x^c u~,~~\kappa (v)= u~,~~\kappa (u)= v~,
\ena
\eq
\epsi (\T{a}{b})=\de^a_b ~,~~\epsi (x^a)=0 ~,~~
\epsi (u)=\epsi (v)=\epsi(I)=1 ~.\label{cfin}
\en
\sk
In the commutative limit $q\rightarrow 1 , r\rightarrow 1$ we
recover the algebra of functions
on $ISO(N)$ and $ISp(N)$ (plus the dilatation $v$ that can be set to
the identity).

{\sl Note 3.1} : as shown in ref. \cite{qpoincare1},
the quantum groups $\ISOqrN$ and $\ISpqrN$
can be derived as the quotients
\eq
\frac{\SOqrNt}{H}~, ~~\frac{\SpqrNt}{H}
\label{quotient}
\en
\noi where $H$ is the Hopf ideal in
 $\SOqrNt$ or $\SpqrNt$ of
all sums of monomials containing at least an element of the kind
$\T{a}{\circ}, \T{\bullet}{b}, \T{\bullet}{\circ}$.
The Hopf structure of the groups in the numerators of
(\ref{quotient})
is naturally inherited by the quotient groups \cite{Sweedler}.

We introduce the following convenient notations:
$\Tc$ stands for
$\T{a}{\circ}$,  $\T{\bullet}{b}$ or  $\T{\bullet}{\circ}$,
$\SqrNt$ stands for either $\SOqrNt$ or
$\SpqrNt$,
and we indicate by
$\Dtwo$, $\epsitwo$ and $\kappatwo$ the corresponding
co-structures.

We denote by $P$ the canonical projection
\eq
P ~:~~  \SqrNt\longrightarrow \SqrNt/{H}
\en
\noi It is a Hopf algebra epimorphism because $H=Ker(P)$ is a Hopf
ideal. Then
any element of ${\SqrNt/H}$ is of the form $P(a)$ and
\eq
P(a)+P(b)\equiv P(a+b) ~;~~ P(a)P(b)\equiv P(ab) ~;~~
\mu P(a)\equiv P(\mu a),~~~\mu \in \mbox{\bf C} \label{isoalgebra}
\en
\eq
\D (P(a))\equiv (P\otimes P)\Dtwo(a) ~;~~ \epsi(P(a))
\equiv\epsitwo(a) ~;~~
\kappa(P(a))\equiv P(\kappatwo(a)) \label{co-iso}
\en
Eq.s (\ref{PRTT11}) - (\ref{PRTT44})  have been obtained in
\cite{qpoincare1}
by taking the $P$
projection of the $RTT$ and $CTT$ relations of
$\SqrNt$, with the notation $u\equiv P(\T{\circ}{\circ}),~v\equiv
P(\T{\bullet}{\bullet}),~z\equiv P(\T{\circ}
{\bullet}), ~x^a \equiv P(\T{a}{\bullet}),~y_a \equiv
P(\T{\circ}{a}),
{}~T^a{}_b\equiv P(T^a{}_b) ~;~~ I\equiv P(I) ~;~ 0\equiv P(0)$,
cf. (\ref{names}).

{\sl Note 3.2 :} From the commutations
(\ref{PRTT22}) - (\ref{PRTT24})
 we see that
one can set $u=I$ only when $q_{a\bullet}=1$ for all $a$.
{}From $q_{a\bullet} = r^2 /q_{\ap\bullet}$, cf. eq. (\ref{qab2}),
this implies also $r=1$.
\sk

{\sl Note 3.3} : eq.s (\ref{PRTT13}) are the multiparametric
(orthogonal or symplectic) quantum plane commutations. They follow
from the   $({}^a{}_{\bullet} {}^b{}_{\bullet})$  $RTT$ components
and (\ref{PRTT44}).
\sk
Finally, the two real forms of $\SqrNt$ mentioned in the previous
section are inherited by $\ISqrN$, with the following
conditions on the parameters:
\sk

$\bullet$~~ $|q_{ab}|=|q_{a\bullet}|=|r|=1$ for $ISO_{q,r}(n,n;\Rb)$,
$ISO_{q,r}(n,n+1;\Rb)$
and $ISp_{q,r}(n;\Rb)$.
\sk
$\bullet$~~For $ISO_{q,r}(n+1,n-1;\Rb)$ : $|r|=1$;
$|q_{a\bullet}|=1$ for $a \not= n, n+1$;
$|q_{ab}|=1$ for $a$
and $b$
both different from $n$ or $n+1$;  $q_{ab}/r \in {\bf R}$
when at least one of the indices $a,b$ is equal
to $n$ or $n+1$;   $q_{a\bullet}/r \in \Rb$ for $a=n$ or $a=n+1$.
\sk
In particular,  the quantum Poincar\'e group $ISO_{q,r}(3,1;\Rb)$
is obtained by setting $|q_{1\bullet}|=|r|=1$,
$q_{2\bullet}/r \in \Rb$, $q_{12}/r  \in \Rb$.

According to {\sl Note 3.1},  a dilatation-free $q$-Poincar\'e group
is found after the further restrictions
$q_{1\bullet}=q_{2\bullet}=r=1$.
The only free parameter remaining is then $q_{12} \in \Rb$.
\sk


\sect{Universal enveloping algebra $U(S_{q,r}(N+2))$}


The universal enveloping algebra of $\SqrNtwo$, i.e. the
algebra of regular functionals \cite{FRT} on $\SqrNtwo$,
is generated
by the functionals $\LLpm $, and the counit  $\epsi$.

The $\LLpm$ linear functionals on $\SqrNtwo$ are defined
by their value on the matrix elements $\T{A}{B}$  :
\eq
\Lpm{A}{B} (\T{C}{D})= \Rpm{AC}{BD}~, ~~~
\Lpm{A}{B} (I)=\de^A_B \label{LonT}
\en
\noi with
\eq
\Rp{AC}{BD} \equiv \R{CA}{DB} ~,~~~
\Rm{AC}{BD} \equiv \Rinv{AC}{BD}~. \label{Rplusminus}
\en
To extend the definition (\ref{LonT})
to the whole algebra $\SqrNtwo$ we set
\eq
\Lpm{A}{B} (ab)=\Lpm{A}{C} (a) \Lpm{C}{B} (b)
{}~~~\forall a,b\in  \SqrNtwo ~.
\label{Lab}
\en

The commutations between $\Lpm{A}{B}$
and $\Lpm{C}{D}$ are given by:
\eq
R_{12} \LLpm_2 \LLpm_1=\LLpm_1 \LLpm_2 R_{12}~,~~~
R_{12} \LLp_2 \LLm_1=\LLm_1 \LLp_2 R_{12}~, \label{RLL}
\en
\noi where as usual the product $\LLpm_2 \LLpm_1$
is the convolution
product $\LLpm_2 \LLpm_1 \equiv (\LLpm_2 \otimes \LLpm_1)\D$.
\sk
{\sl Note 4.1 :} $L^+$ is upper
triangular and $L^-$ is lower triangular. Proof:
apply $L^+ $ and  $L^-$ to
the $T$ elements and use the upper and lower
triangularity of $R^+$ and $R^-$, respectively.
\sk
The $\Lpm{A}{B}$ elements satisfy orthogonality conditions
analogous to (\ref{CTTbig}):
\eq
C^{AB} \Lpm{C}{B} \Lpm{D}{A} = C^{DC} \epsi,~~~
C_{AB} \Lpm{B}{C} \Lpm{A}{D} = C_{DC} \epsi \label{CLL}
\en
\noi as can be verified by applying them to the $q$-group
generators, and using (\ref{crc}).
They provide a quantum inverse for $\Lpm{A}{B}$:
\eq
(\Lpm{A}{B})^{-1}=C^{DA} \Lpm{C}{D} C_{BC} \label{Linverse}
\en
\sk
The co-structures of the algebra generated by the
functionals $L^{\pm}$ and $\epsi$
are defined by:
\eq
\Dp(\Lpm{A}{B})(a \otimes b) \equiv \Lpm{A}{B}
(ab)=\Lpm{A}{G}(a) \Lpm{G}{B} (b)~,
\en
\eq
\ep (\Lpm{A}{B})\equiv \Lpm{A}{B} (I)~~;~~~
\kp  (\Lpm{A}{B})(a)\equiv \Lpm{A}{B} (\kappa (a)) =
(\Lpm{A}{B})^{-1} (a)~,
\en
\noi so that
\eqa
& &\Dp (\Lpm{A}{B})=\Lpm{A}{G} \otimes \Lpm{G}{B}~,\label{copLpm}\\
& &\ep (\Lpm{A}{B})=\de^A_B ~~;~~~
\kp (\Lpm{A}{B})= (\Lpm{A}{B})^{-1}= C^{DA} \Lpm{C}{D} C_{BC}~.
\label{coiLpm}
\ena
\sk
The *-conjugation on $\SqrNtwo$ induces a *-conjugation
on $\USqrNtwo$ in two possible ways
(we denote them as $*$ and $\sharp$-conjugations):
\eq
\phi^* (a)\equiv {\overline {\phi (\kappa  (a^*))}}
\label{starconjugation}~~~;~~~~
\phi^{\sharp} (a)\equiv {\overline {\phi (\kappa^{-1} (a^*))}}
\en
\noi where $\phi \in \USqrNtwo$, $a \in \SqrNtwo$, and the overline
denotes
 the usual complex conjugation.
Both $*$ and $\sharp$ can be shown
to satisfy all the properties of Hopf algebra involutions.
It is not difficult
to determine their action
on the basis elements $\Lpm{A}{B}$.
The two $\SqrNtwo$ $*$-conjugations  of
the previous section induce
respectively the following conjugations on the
$\Lpm{A}{B}$:
\eqa
& &~~~~(\Lpm{A}{B})^*=\Lpm{A}{B} \label{conjL1}\\
& & ~~~(\Lpm{A}{B})^*=\Dcal{A}{C} \Lpm{C}{D} \Dcal{D}{B}
\label{conjL2}
\ena
\noi To find $(\Lpm{A}{B})^{\sharp}$ one uses the general formula
$(\phi)^{\sharp}=\kpsq~[ (\phi)^*]$, deducible from the compatibility
of both conjugations with the antipode: ${\kp}^{-1} (\phi^*)=[\kp
(\phi)]^*,~
 {\kp}^{-1} (\phi^\sharp)=[\kp (\phi)]^\sharp$.


\sect{Bicovariant calculus on $\SqrNtwo$}

The bicovariant differential calculus on the uniparametric
$q$-groups of the $A,B,C,D$ series can be formulated in terms
of the corresponding $R$-matrix, or equivalently in terms of
the $\LLpm$ functionals.  This holds also for the multiparametric
case.  In fact all formulas are the same, modulo substituting
the $q$ parameter with $r$  when
it  appears explicitly (typically as ${1 \over {q-q^{-1}}}$).
\sk
We briefly recall how to construct a bicovariant
calculus.  The general procedure can be found in ref.
\cite{Jurco,Schmudgen}, or, in the notations we adopt
here, in ref. \cite{Aschieri1}.
It realizes the axiomatic construction of ref. \cite{Wor}.
\sk

As in the uniparametric case \cite{Jurco}, the functionals
\eq
\ff{A_1}{A_2B_1}{B_2} \equiv \kp (\Lp{B_1}{A_1}) \Lm{A_2}{B_2}.
\label{defff}
\en
and the elements of $A=\SqrNtwo$:
\eq
\MM{B_1}{B_2A_1}{A_2} \equiv \T{B_1}{A_1} \kappa (\T{A_2}{B_2}).
\label{defMM}
\en
satisfy the following relations, called bicovariant bimodule
conditions,
where for simplicity
we use the adjoint indices $ i,j,k,...$
with ${}^i={}_A^{~B},~{}_i={}^A_{~B}~:$
\eqa
& & \D '(\f{i}{j}) = \f{i}{k} \otimes \f{k}{j}  \label{propf1}~~;~~~
\epsi '(\f{i}{j}) = \del^i_j~,  \\
& & \Delta (\M{j}{i}) = \M{j}{l} \otimes \M{l}{i} \label{copM}~~;~~~
\epsi (\M{j}{i}) = \delta^i_j~, \\
& & \M{i}{j} (a * \f{i}{k})=(\f{j}{i} * a) \M{k}{i} \label{propM}
\ena
\noi The star product between a functional on $A$ and an
element of $A$ is defined as:
\eq
 \chi * a \equiv (id \otimes \chi) \D (a),~~~
a * \chi \equiv (\chi \otimes id) \D (a),~~~
a \in A, ~\chi \in A'
\en
Relation  (\ref{propM}) is easily checked for
$a=\T{A}{B}$ since in
this case it is implied by the $RTT$ relations;
it holds for a generic $a$ because of  property
(\ref{propf1}).
\sk
The space of {\bf  quantum one-forms} is defined as a right
$A$-module $\Ga$
freely generated  by the  {\sl symbols} $\ome{A_1}{A_2}$:
\eq
\Ga\equiv \{a^{A_1}_{~A_2}\ome{A_1}{A_2}\}~ ,~~~a^{A_1}_{~A_2} \in
A\label{symbols}
\en
{\sl Theorem 5.1} (due to Woronowicz: see Theorem 2.5  in
the last ref. of \cite{Wor},
p. 143): because of the
properties  (\ref{propf1}), $\Ga$
becomes a bimodule over $A$ with the following
right multiplication:
\eq
\ome{A_1}{A_2} a=(\ff{A_1}{A_2B_1}{B_2} * a)~\ome{B_1}{B_2},
\label{omea}
\en
in particular:
\eq
\ome{A_1}{A_2} \T{R}{S}= \Rinv{TB_1}{CA_1}
\Rinv{A_2C}{B_2S} \T{R}{T}\ome{B_1}{B_2} \label{commomT}
\en
Moreover, because of properties (\ref{copM}),
 we can  define a left and a right action of
$A$ on $\Ga$:
\eqa
{}\!\!\!\!\!\!
& & \DL:\Ga\rightarrow A\otimes
\Ga~~;~~~\DL(a \ome{A_1}{A_2} b)
\equiv \D(a)\,(I\otimes \ome{A_1}{A_2})\,\D(b)~,
\label{deltaLomega} \\
{}\!\!\!\!\!\!& & \DR:\Ga\rightarrow  \Ga\otimes A
{}~~;~~~\DR(a \ome{A_1}{A_2} b)
\equiv \D(a)\,(\ome{B_1}{B_2}\otimes
\MM{B_1}{B_2A_1}{A_2})\,\D(b)~.
\label{deltaRomega}
\ena
These actions commute, i.e. $(id \otimes \DR) \DL = (\DL \otimes id)
\DR$ because of (\ref{propM}),
and give a bicovariant bimodule structure to $\Ga$.
\sk
The {\bf exterior derivative}
$d ~:~ A\longrightarrow \Ga$ can be defined via the element
$\tau\equiv \sum_A \ome{A}{A} \in \Ga$. This element
is
easily shown to be
left and right-invariant:
\eq
\DL(\tau)=I\otimes\tau~~;~~~
\DR(\tau)=\tau\otimes I
\en
and defines the derivative  $d$  by
\eq
da=\rinv [\tau a - a \tau]. \label{defd1}
\en
\noi The factor $\rinv$ is necessary for a correct
classical limit $\rone$.
It is immediate to prove the Leibniz rule
\eq
d(ab)=(da)b+a(db),~~\forall a,b\in A ~. \label{Leibniz}
\en
Another expression for the derivative is given by:
\eq
da = ( \cchi{A_1}{A_2} * a) ~\ome{A_1}{A_2} \label{defd2}
\en
where the linearly independent elements
\eq
\cchi{A}{B} = \rinv [\ff{C}{CA}{B}-\de^A_B \epsi]
\label{defchi}
\en
are the tangent vectors such that the left-invariant vector fields  
$\chi^A{}_B
*$
are dual to the left-invariant
one-forms $\ome{A_1}{A_2}$.
The equivalence of  (\ref{defd1}) and (\ref{defd2}) can be
shown by using the rule (\ref{omea}) for $\tau a$ in
the right-hand side of (\ref{defd1}).

\noi Using (\ref{defd2}) we compute the exterior derivative
on the basis elements of $\SqrNtwo$:
\eq
d ~\T{A}{B}=\rinv [\Rinv{CR}{ET} \Rinv{TE}{SB}
\T{A}{C}-\de^R_S \T{A}{B}]
{}~\ome{R}{S} \equiv\T{A}{C}\X{CR}{BS}
 \ome{R}{S}\label{dTAB}
\en
where we have
\eq
\X{A_1B_1}{A_2B_2}\equiv \rinv [\Rinv{A_1 B_1}{ET} \Rinv{TE}{B_2A_2}-
\de^{B_1}_{B_2} \de^{A_1}_{A_2} ]
=zK^{A_1B_1}_{~~A_2B_2} -\Rhatinv{A_1B_1}{A_2B_2}
\label{X}
\en
with  $z \equiv \epsilon r^{N-\epsilon},
K^{A_1B_1}_{~~A_2B_2}=C^{A_1B_1}C_{A_2B_2} ~.$ [From
(\ref{extrarelation}), the second equality in (\ref{X}) is easily  
proven.]
Every element  $\rho$ of $\Gamma$, which by definition is
written in a unique way as
$\rho=a^{A_1}{}_{A_2}\ome{A_1}{A_2}$,  can also be written as
\eq
\rho=\sum_k a_k db_k  \label{propd}
\en
\noi for some $a_k,b_k$ belonging to $A$. This can
be proven directly by inverting
the relation (\ref{dTAB}). The
result  is an expression of
 the $\om$ in terms of a linear
combination of $\kappa (T) dT$, as in the classical case:
\eq
\ome{A_1}{A_2}=Y_{A_1B_1}^{~~~A_2B_2} \kappa(\T{B_1}{C}) d\T{C}{B_2}
\label{omAB}
\en
\noi where $Y$ satisfies
$
\X{A_1B_1}{A_2B_2}Y_{B_1C_1}^{~~~B_2C_2}=\de^{A_1}_{C_1}\de^{C_2}_{A_2 
} ~,~
Y_{A_1B_1}^{~~~A_2B_2} \X{B_1  
C_1}{B_2C_2}=\de^{C_1}_{A_1}\de^{A_2}_{C_2}
$
and is given explicitly by
\eq
Y_{A_1B_1}^{~~~A_2B_2}=
\alpha[(z-\la) C_{A_1B_1}C^{A_2B_2} + C_{A_1D}
\R{DA_2}{CB_1} C^{CB_2}
 - {\la\over z(z-z^{-1})}
D^{A_2}_{~A_1} (D^{-1})^{B_2}_{~B_1}] \label{Y}
\en
\noi with $\alpha={1 \over z(z-z^{-1}-\la)}$ and
$D^{E}_{~C} \equiv C^{EF} C_{CF}$.
The $r=1$ limit of  (\ref{dTAB}) is discussed in the next section.
\sk
Due to the bi-invariance of $\tau$ the derivative operator $d$ is
compatible with the actions $\DL$ and $\DR$:
\eq
\DL (adb)=\D(a)(id\otimes d)\D(b)~~,~~~\DR (adb)=\D(a)(d\otimes
id)\D(b)~,\label{dbicov}
\en
these two properties express the fact that $d$
commutes with the left
and right action of the quantum group, as in the classical case.
\sk
\noi {\sl Remark 5.1:} The properties (\ref{Leibniz}),
(\ref{propd})  and (\ref{dbicov}) of
the exterior derivative (\ref{defd2})
realize the axioms of a first-order bicovariant differential
calculus \cite{Wor}.
\sk
The {\bf tensor product} between elements $\rho,\rhop \in \Ga$
is defined to
have the properties $\rho a\otimes \rhop=\rho
\otimes a \rhop$, $a(\rho
\otimes \rhop)=(a\rho) \otimes \rhop$ and
$(\rho \otimes \rhop)a=\rho
\otimes (\rhop a)$. Left and right actions on
$\Ga \otimes \Ga$ are
defined by:
\eqa
{}\!\!\!\!\!\!\!
&&\DL (\rho \otimes \rhop)\equiv \rho_1 \rhop_1
\otimes \rho_2 \otimes
\rhop_2,~~~\DL: \Ga \otimes \Ga \rightarrow A
\otimes\Ga\otimes\Ga
\label{DLGaGa}\\
{}\!\!\!\!\!\!\!
&&\DR (\rho \otimes \rhop)\equiv \rho_1 \otimes \rhop_1
 \otimes \rho_2
\rhop_2,~~~\DR: \Ga \otimes \Ga \rightarrow
\Ga\otimes\Ga\otimes A
\label{DRGaGa}
\ena
\noi where  $\rho_1$, $\rho_2$, etc., are defined by:
$$
\DL (\rho) = \rho_1 \otimes \rho_2,~~~\rho_1\in A,~\rho_2\in
\Ga~~;~~~
\DR (\rho) = \rho_1 \otimes \rho_2, ~\rho_1\in \Ga,~\rho_2\in A~.$$
The extension to $\Ga^{\otimes n}$ is straightforward.
\sk
The  {\bf exterior product} of one-forms is consistently defined as:
\eq
\ome{A_1}{A_2} \we \ome{D_1}{D_2}
\equiv \ome{A_1}{A_2} \otimes \ome{D_1}{D_2}
- \RRhat{A_1}{A_2}{D_1}{D_2}{C_1}{C_2}{B_1}{B_2}
\ome{C_1}{C_2} \otimes \ome{B_1}{B_2} \label{exteriorproduct}
\en
\noi where the $\Lambda$ tensor is given by:
\eq
\begin{array}{rl}
\LL{A_1}{A_2}{D_1}{D_2}{C_1}{C_2}{B_1}{B_2}
\equiv & \ff{A_1}{A_2B_1}{B_2} (\MM{C_1}{C_2D_1}{D_2}) = \\
=&
d^{F_2} d^{-1}_{C_2} \R{F_2B_1}{C_2G_1} \Rinv{C_1G_1}{E_1A_1}
    \Rinv{A_2E_1}{G_2D_1} \R{G_2D_2}{B_2F_2} \label{Lambda}
\end{array}
\en
This matrix satisfies the characteristic equation:
\eq
\begin{array}{rcl}
& & (\La+r^2 I)~(\La+r^{-2}I)~(\La + \epsilon r^{\epsilon + 1-N} I)
(\La + \epsilon r^{-\epsilon - 1+N} I) \times\\
& &~~~~~~~(\La - \epsilon r^{-\epsilon + 1+N} I)
{}~(\La - \epsilon r^{\epsilon - 1-N} I)~
(\La-I)=0
\end{array} \label{Laeigen}
\en
\noi due to the characteristic equation  (\ref{cubic}).
For simplicity
we will at times use the adjoint indices $i,j,k,...$
with ${}^i={}_A^{~B},~{}_i={}^A_{~B}$.
Define
\eq
\PIJ{a_1}{a_2}{d_1}{d_2}{c_1}{c_2}{b_1}{b_2}
\equiv  d^{f_2} d^{-1}_{c_2}
\Rhat{b_1f_2}{c_2g_1} \PI{c_1g_1}{a_1e_1}
\Rhatinv{a_2e_1}{d_1g_2} \PJ{d_2g_2}{b_2f_2} \label{PIJ}
\en
\noi where  $P_I=P_S,P_A,P_0$ are given in
(\ref{projBCD}). The $(P_I,P_J)$ are themselves
 projectors, i.e.:
\eq
(P_I,P_J) (P_K,P_L) =\de_{IK} \de_{JL} (P_I,P_J) \label{projIJ}
\en
Moreover
\eq
(I,I)=I \label{IIeqI}
\en
{}From (\ref{exteriorproduct}) we find
\eq
\om^i \we \om^j = - \Z{ij}{kl} \om^k \we \om^l \label{commom}
\en
\noi with
 \eq
Z=(P_S,P_S)+(P_A,P_A)+(P_0,P_0)-I
\en
\noi see ref. \cite{qconstants}.
The inverse of $\Lambda$ always
exists, and is given by
\eq
\begin{array}{rl}
\!\!\!\RRhatinv{A_1}{A_2}{D_1}{D_2}{B_1}{B_2}{C_1}{C_2}=&
\ff{D_1}{D_2B_1}{B_2}
(\T{A_2}{C_2} \km (\T{C_1}{A_1})) = \\
=&\R{F_1B_1}{A_1G_1} \Rinv{A_2G_1}{E_2D_1} \Rinv{D_2E_2}{G_2
C_2} \R{G_2C_1}{B_2F_1} (d^{-1})^{C_1} d_{F_1}
\end{array}
\label{Lambdainv}
\en
\noi Note that for $r=1$, $\La^2 = I$ and $(\La +I)(\La-I)=0$
replaces the
seventh-order spectral equation (\ref{Laeigen}). In this special
case, one
finds
the simple formula:
\eq
\om^i \we \om^j= - \L{ij}{kl} \om^k \we \om^l
\label{omcomrone}~~~~\mbox{ i.e.}
{}~~Z=\Lambda~.
\en
\sk
Using the exterior product we can define the {\bf exterior
differential on $\Ga$} :
\eq
d~:~\Gamma \rightarrow \Gamma \we \Gamma
{}~~~;~~~~d (a_k db_k) = da_k \we db_k \label{defdonga}
\en
\noi which can easily be extended to
$\Gamma^{\we n}$ ($d: \Gamma^{\we n}
\rightarrow \Gamma^{\we (n+1)}$, $\Gamma^{\we n}$ being
defined as in the
classical case but with the quantum permutation
(braid) operator $\La$ \cite{Wor}). The definition (\ref{defdonga})
is equivalent to the following :
\eq
d\theta = \rinv [\tau \we \theta - (-1)^k \theta \we \tau],
\label{defdgen}
\en
\noi where $\theta \in \Ga^{\we k}$. The exterior differential
 has the following properties:
\eq
d(\theta \we \thetap)=d\theta \we \thetap +
(-1)^k \theta \we d\thetap~~~;~~~~
d(d\theta)=0~,\label{propd2}
\en
\eq
\DL (\theta d\theta')=\DL(\theta)(id\otimes
d)\DL(\theta')\label{propd3}~~~;~~~~
\DR (\theta d\theta')=\DR(\theta)(d\otimes
id)\DR(\theta'),\label{propd4}
\en
\noi where $\theta \in \Ga^{\we k}$, $\thetap \in \Ga^{\we n}$.
\sk
The {\bf \mbox{\boldmath $q$}
-Cartan-Maurer equations} are found by
using
(\ref{defdgen}) in computing $d\ome{C_1}{C_2}$:
\eq
d\ome{C_1}{C_2}= \rinv (\ome{B}{B} \we \ome{C_1}{C_2} +
 \ome{C_1}{C_2} \we
\ome{B}{B}) \equiv
-\onehalf \cc{A_1}{A_2}{B_1}{B_2}{C_1}{C_2}
{}~\ome{A_1}{A_2} \we \ome{B_1}{B_2} \label{CartanMaurer}
\en
\noi with:
\eq
\cc{A_1}{A_2}{B_1}{B_2}{C_1}{C_2}= -{2\over (r-r^{-1})}
[ \ZZ{B}{B}{C_1}{C_2}{A_1}{A_2}{B_1}{B_2} + \de^{A_1}_{C_1}
\de^{C_2}_{A_2} \de^{B_1}_{B_2} ]\label{cc}
\en
To derive this formula we have used the flip operator $Z$
on $\ome{B}{B} \we \ome{C_1}{C_2}$.
\sk

Finally, we recall that the $\chi$ operators close
on the {\bf $q$-Lie algebra} :
\eq
\chi_i \chi_j - \L{kl}{ij} \chi_k \chi_l = \C{ij}{k} \chi_k
\label{bico1}
\en
\noi where the $q$-structure constants are given by
\eq
{\C{jk}{i}=\chi_k(\M{j}{i})} ~~\mbox{ explicitly : }~~
{\CC{A_1}{A_2}{B_1}{B_2}{C_1}{C_2}} =\rinv [- \de^{B_1}_{B_2}
\de^{A_1}_{C_1} \de^{C_2}_{A_2} +
\LL{B}{B}{C_1}{C_2}{A_1}{A_2}{B_1}{B_2}]. \label{CC}
\en
The $C$ structure constants
appearing in the Cartan-Maurer equations are in general related
to the $\Cb$ constants of the $q$-Lie algebra \cite{Aschieri1}:
\eq
\C{jk}{i}=\onehalf [\c{jk}{i}-\L{rs}{jk} \c{rs}{i}]~.
\en
\noi In the particular case $\Lambda^2=I$ (i.e. for $r=1$) it is not
difficult to
see that in fact
$C= \Cb$, and that
the $q$-structure constants are $\Lambda$-antisymmetric:
\eq
\C{jk}{i}= - \L{rs}{jk} \C{rs}{i}. \label{antisymC}
\en

The $\chi$ and $f$ operators close on the algebra
(\ref{bico1}) and
\eq
\L{nm}{ij} \f{i}{p} \f{j}{q} = \f{n}{i} \f{m}{j} \L{ij}{pq}
\label{bico2}
\en
\eq
\C{mn}{i} \f{m}{j} \f{n}{k} + \f{i}{j} \chi_k= \L{pq}{jk} \chi_p
\f{i}{q} + \C{jk}{l} \f{i}{l} \label{bico3}
\en
\eq
\chi_k  \f{n}{l}=\L{ij}{kl} \f{n}{i} \chi_j~,  \label{bico4}
\en
This algebra is {\sl sufficient} to define
a bicovariant differential calculus on $A$
(see e.g. \cite{Bernard}), and will be called
``bicovariant algebra" in the sequel.
By applying the relations defining
the bicovariant algebra
to the element $\M{r}{s}$
we can express them
in the adjoint representation:
\eqa
& & \C{ri}{n} \C{nj}{s}-\L{kl}{ij} \C{rk}{n} \C{nl}{s} =
\C{ij}{k} \C{rk}{s}
{}~~\mbox{({\sl q}-Jacobi identities)} \label{bicov1}\\
& & \L{nm}{ij} \L{ik}{rp} \L{js}{kq}=\L{nk}{ri} \L{ms}{kj}
\L{ij}{pq}~~~~~~~~~\mbox{(Yang--Baxter)} \label{bicov2}\\
& & \C{mn}{i} \L{ml}{rj} \L{ns}{lk} + \L{il}{rj} \C{lk}{s} =
\L{pq}{jk} \L{is}{lq} \C{rp}{l} + \C{jk}{m} \L{is}{rm}
\label{bicov3}\\
& & \C{rk}{m} \L{ns}{ml} = \L{ij}{kl} \L{nm}{ri} \C{mj}{s}
\label{bicov4}
\ena
\sk
Using the definitions (\ref{defchi}) and (\ref{defff}) it is
not difficult to find the co-structures
on the functionals $\chi$ and $f$:
\eq
\begin{array}{lcl}
\Dp (\chi_i)=\chi_j
\otimes \f{j}{i} + \epsi \otimes \chi_i
&~;~~ &\Dp (\f{i}{j})=\f{i}{k} \otimes \f{k}{j}~, \\
\ep(\chi_i)=0
&~;~~ & \ep (\f{i}{j}) = \del^i_j~, \\
\kp(\chi_i)=-\chi_j \kp(\f{j}{i})
&~;~~ & \kp (\f{k}{j}) \f{j}{i}= \de^k_i \epsi =
\f{k}{j} \kp (\f{j}{i})~.
\end{array}\label{copf}
\en
\noi Note
that in the $r,q \rightarrow 1$ limit $\f{i}{j}
\rightarrow \de^i_j \epsi$, i.e. $\f{i}{j}$ becomes
proportional to the
identity functional and formula (\ref{omea}),
becomes trivial, e.g. $\om^i a = a\om^i$ [use $\epsi * a=a$].
\sk
\noi{\sl Note 5.1:} The formulae characterizing the bicovariant  
calculus have
been written in the
 basis $\{\cchi{A}{B}\},~ \{\ome{C}{D}\}$  because of the  
particularly simple
expression of the
$\ff{A}{BC}{D}$ and $\cchi{A}{B}$ functionals
in terms of $L^{\pm A}{}_B$, see (\ref{defff}) and (\ref{defchi}).
Obviously the calculus is independent from the basis chosen. If we  
consider the
linear transformation
\[
\om^i\rightarrow \om '^i=X^i{}_j\om^j
\]
(where we use adjoint indices ${~}^i={}_{A_1}{}^{A_2} ,
{}_j={}^{B_1}{}_{B_2}$), from the exterior differential
\eq
da=(\chi_i*a)\om^i=(\chi'_i*a)\om '^i
\en
we find
\[
\chi_i\rightarrow\chi'_i=\chi_j(X^{-1}){}^j{}_i~,
\]
and from the coproduct rule (\ref{copf}) of the $\chi_i$ we find
$f^{i}{}_{j}\rightarrow f'^i{}_j=X^i{}_lf^l{}_m(X^{-1}){}^m{}_j ;$  
while from
(\ref{deltaRomega})
we have $M_i{}^j\rightarrow M'_i{}^j=(X^{-1}){}^l{}_iM_l{}^mX^j{}_m  
.$

A useful change of basis is obtained via the following  
transformation:
\eq
\begin{array}{rcl}
&  
&\ome{A_1}{A_2}\rightarrow\vt^{A_1}_{~A_2}=\X{A_1B_1}{A_2B_2}\ome{B_1} 
{B_2}\\
& &\cchi{A_1}{A_2}\rightarrow
\psi_{A_1}^{~A_2}=\cchi{B_1}{B_2}Y_{B_1A_1}^{~~B_2A_2}
\end{array}
\en
where $X$ and its (second) inverse $Y$ are defined in (\ref{X}) and  
(\ref{Y}).
Using (\ref{omAB}) it is immediate to see that
\eq
\vt^{A_1}_{~A_2}=\kappa(\T{A_1}{C})d\T{C}{A_2}~.
\en
We also have\footnote{We recall from \cite{Wor} (\cite{PaoloPeter})  
that the
quantum group elements (coordinates)
$x^j$ such that
$$
x^j\in Ker\,\varepsilon  ~~\mbox{ and }~~\chi_i(x^j)=\de_i^j
$$
are uniquely defined by these two conditions.
Notice, by the way, that
$f^i{}_j(a)=\chi_j(x^ia)~.$
} :
\eq
\psi_{A_1}^{~A_2}({T}^{B_1}_{~B_2})=
\psi_{A_1}^{~A_2}({\widetilde  
T}^{B_1}_{~B_2})=\de_{A_1}^{B_1}\de^{A_2}_{B_2}~~
\mbox{ where }~~ {\widetilde T}^{B_1}_{~B_2}\equiv
\T{B_1}{B_2}-\de^{B_1}_{B_2}I~.\label{dualchix}
\en
Formula (\ref{dualchix}) follows from $\psi_{A_1}^{~A_2}(I)=0$ and:
\eq
\begin{array}{rcl}
\vt^{A_1}_{~A_2}&=&\kappa(\T{A_1}{C})d\T{C}{A_2}=\kappa(\T{A_1}{C})
(\psi_{B_1}^{~B_2}*\T{C}{A_2})\vt^{B_1}_{~B_2}\\
&=&\kappa(\T{A_1}{C})\T{C}{D}
\psi_{B_1}^{~B_2}(\T{D}{A_2})\vt^{B_1}_{~B_2}=\psi_{B_1}^{~B_2}(\T{A_1 
}{A_2})
\vt^{B_1}_{~B_2}~.
\end{array}
\en
The analogue of the coordinates ${\widetilde T}^{B_1}_{~B_2}$ in the  
old basis
is given by
\eq
x_{B_1}^{~B_2}\equiv Y_{B_1C_1}^{~~B_2C_2}{\widetilde
T}^{C_1}_{~C_2}~~,~~~\cchi{A_1}{A_2}
(x_{B_1}^{~B_2})=\de^{A_1}_{B_1}\de^{B_2}_{A_2}~.
\en
\sk

Compatibility of the {\bf conjugation} defined in
(\ref{starconjugation})
with the differential calculus requires
$(\chi_i)^*$ to be a linear combination
of $(\kp)^{-2} (\chi_i)$, or
$(\chi_i)^{\sharp}$ to be a linear combination
of the $\chi_i$.
This follows from Theorem 1.10 (Woronowicz)  of  last ref. in
\cite{Wor},
and from equations (\ref{starconjugation})
with $\phi=\chi$.

{}From the definitions (\ref{defff}),
(\ref{defchi}), it is straightforward to find
how the $*$ and $\sharp$ conjugations act
on the tangent vectors $\chi$. Both the conjugations (\ref{conjL1})
and (\ref{conjL2})
satisfy the criteria given above for its compatibility
with the differential calculus. Indeed the conjugation
(\ref{conjL2}) yields [use (\ref{defchi}),  (\ref{defff}),   
(\ref{coiLpm}),
(\ref{RLL})]:
\eq
(\cchi{A}{B})^*=-\epsilon r^{1-N} \cchi{C}{D} \Dcal{F}{~B}
 \Dcal{A}{~G} \R{DG}{FE} D^E_{~C} \label{conjchiAB}
\en
\noi with $D^E_{~C} \equiv C^{EF} C_{CF}$. To find
$\chi^*$ corresponding to (\ref{conjL1}) just take
$\Dcal{A}{B}=\de^A_B$ in (\ref{conjchiAB}). The criterion
given above for the compatibility with the differential calculus is
fulfilled since $(\kp)^{-2} (\chi_i)$ is a linear combination
of $\chi_i$:
\eq
\kpsq (\cchi{A}{B})= (D^{-1})^{A}_{~C} \cchi{C}{D} D^{D}_{~B}
\en
\noi as can be seen from  (\ref{defchi}) and $\kpsq (\Lpm{A}{B})=
 (D^{-1})^{A}_{~C} \Lpm{C}{D} D^{D}_{~B} $, cf. (\ref{coiLpm}).

Using the inversion formulae (\ref{omAB}) one finds
the induced conjugation on the
left-invariant one-forms:
\eq
(\ome{A_1}{A_2})^*= -\Dcal{F}{B_2} \Dcal{B_1}{G} C_{FA} C^{MG}
\overline{ (Y_{A_1B_1}^{~~A_2B_2})} \X{AC_1}{MC_2} \ome{C_1}{C_2}~.
\label{conjomAB}
 \en
\sk


\sect{Differential calculus on $\SqroNt$}

As discussed in Section 3,  we have obtained the quantum
inhomogeneous groups $\ISqrN$ via the projection
\eq
P~:~~~\SqrNtwo -\!\!\!-\!\!\!\!\!\longrightarrow {\SqrNtwo \over
H}=\ISqrN
\label{Pproj}
\en
with $H$=Hopf ideal in $\SqrNtwo$ defined in Section 3.  As a
consequence, the universal enveloping algebra
$\UISqrN$ is a Hopf subalgebra of $\USqrNtwo$
\cite{Sweedler,UEA},
and contains all the functionals that annihilate $H=Ker(P)$.

Let us consider now the $\chi$ functionals in the differential
calculus
on $\SqrNtwo$. Decomposing the indices we find:
\eqa
& & \cchi{a}{b}=\rinv [\ff{c}{c a}{b}-\de^a_b \epsi]~~~~~~~+\rinv
\ff{\bu}{\bu
a}{b} \label{chiN2}\\
& & \cchi{a}{\ci}=\rinv \ff{c}{c a}{\ci} ~~~~~~~~~ ~~~~~~~ +\rinv
\ff{\bu}{\bu
a}{\ci} \\
& &\cchi{\ci}{b}=~~~~~~~~~~~~~~~~~~~~~~~~~~~~~~~~~~ + \rinv
 [\ff{c}{c \ci}{b}+\ff{\bu}{\bu \ci}{b}]  \\
& & \cchi{a}{\bu}=~~~~~~~~~~~~~~~~~~~~~~~~~~~~~~~~~~+ \rinv
\ff{\bu}{\bu a}{\bu}\\
& & \cchi{\bu}{b}=\rinv \ff{\bu}{\bu\bu}{b} ~~~~~~~~~~~~~\\
& &\cchi{\ci}{\ci}=\rinv [\ff{\ci}{\ci\ci}{\ci}-\epsi]~~~~~~~~~~ +
\rinv
[\ff{c}{c\ci}{\ci}
+\ff{\bu}{\bu \ci}{\ci} ] \\
& &\cchi{\ci}{\bu}=~~~~~~~~~~~~~~~~~~~~~~~~~~~~~~~~~~ + \rinv
\ff{\bu}{\bu \ci}{\bu}\\
& & \cchi{\bu}{\ci}=\rinv \ff{\bu}{\bu \bu}{\ci}~~~~~~~~~~~~~~~\\
& & \cchi{\bu}{\bu} = \rinv [\ff{\bu}{\bu\bu}{\bu} -
\epsi]~~~~~\label{chiN2end}\\
& &~~~~~~~~\underbrace{~~~~~~~~~~~~~~~~~~~~~~~~~~~~~~~~~}_
{\hbox{terms annihilating $H$}} \nonumber
\ena
\noi where we have indicated the terms that do and do not
annihilate the Hopf ideal
$H$. We see that only the functionals
$\cchi{\bu}{b}$,  $\cchi{\bu}{\ci} $ and $\cchi{\bu}{\bu}$ do
annihilate
$H$, and therefore belong to $\UISqrN$. The resulting bicovariant
differential calculus \cite{UEA} contains dilatations and
translations,
but does not contain the tangent vectors of $\SqrN$, i.e.
the functionals $\cchi{a}{b}$.  Indeed these contain
 $\ff{\bu}{\bu a}{b}$, in general not vanishing on $H$. We can,
however, try to find restrictions on the parameters
$q,r$ such that  $\ff{\bu}{\bu a}{b} (H)=0$. As we will see, this
happens
 for $r=1$. For this reason we consider in the following
the particular multiparametric deformations called ``minimal
deformations" or twistings, corresponding to $r=1$.
\sk
We first examine what happens to the
bicovariant calculus on $\SqrNtwo$ in the $r=1$
limit\footnote{By $\lim_{\rone}a$,
where the generic element $a\in\SqrNtwo$  is a polynomial in the
matrix
elements $\T{A}{B}$ with complex coefficients $f(r)$ depending on  
$r$,
we
understand the
element of $S_{q,r=1}(N+2)$ with coefficients given by
$\lim_{\rone}f(r)$. The
expression
$\lim_{\rone}\phi=\varphi$, where $\phi\in \USqrNtwo$  and  
$\varphi\in
U(S_{q,r=1}(N+2))$ means that
$\lim_{\rone}\phi(a)=\varphi(\lim_{\rone}a)$ for any $a\in \SqrNtwo$  
such
that
$\lim_{\rone}a$ exists.
Finally, the left invariant one-forms $\omega^i$ are symbols, and
therefore
$\lim_{\rone}a_i\omega^i\equiv(\lim_{\rone}a_i)\omega^i$ [see
(\ref{symbols})].}.
The $R$ matrix is given by, cf. (\ref{Rnonzero}):
\eqa
& & \R{AB}{AB}=q^{-1}_{AB} + O(\la) \\
& & \R{AB}{BA}=\la~~~~~~~~~~~~~~~~~~~~~~~~
\mbox{\footnotesize $A>B, A' \not= B$}\\
& & \R{AA'}{A'A}=\la~ (1-\epsilon r^{\rho_A - \rho_{A'}})~~~~~
\mbox{\footnotesize  $A>A'$} \\
& & \R{AA'}{BB'}=-\la \epsilon_A \epsilon_B + O(\la^2)~~~~
\mbox{\footnotesize  $A>B, A' \not= B$}
\ena
\noi where $O(\la^n)$ indicates an infinitesimal of order
$\geq\la^n$;
the $q_{AB}$ parameters satisfy:
\eq
q_{AB}=q^{-1}_{AB'} = q^{-1}_{A'B} = q^{-1}_{BA}~~;~~~~ q_{AA} =
q_{AA'} = 1
\en
\noi up to order $O(\la)$.  Note that the components $\R{AA'}{A'A}$
are of order $O(\la^2)$ for the orthogonal case ($\epsilon=1$) and
of order $O(\la)$ for the symplectic case ($\epsilon=-1$).
The $RTT$ relations simply become:
\eq
\T{B_1}{A_1} \T{B_2}{A_2}={q_{B_1B_2} \over q_{A_1A_2}}\T{B_2}{A_2}
\T{B_1}{A_1}~.
\en
For $r=1$ the metric is $C_{AB}=\epsilon_A\delta_{AB'}$
and therefore we have
$C_{AB}=\epsilon C_{BA}$.
Using the
definition
(\ref{LonT}), it is easy to see that
\eqa
{}\!\!\!\!\!{}\!\!\!\!\!&{}\!\!\!\!\!{}\!\!\!\!\! &
\Lpm{A}{A} (\T{C}{D})=\de^C_D q_{AC} +O(\la)\\
{}\!\!\!\!\!{}\!\!\!\!\!&{}\!\!\!\!\!{}\!\!\!\!\! &
\Lpm{A}{B} (\T{B}{A})=\pm \la {}~~~~~~~\mbox{\footnotesize $A\not=B,
A'\not= B$;~$A<B$ for $L^+$, $A>B$ for $L^-$}\label{pmlanda1}\\
{}\!\!\!\!\!{}\!\!\!\!\!&{}\!\!\!\!\!{}\!\!\!\!\! &
\Lpm{A}{A'} (\T{A'}{A})=\pm \la~[1-\epsilon r^{\pm
(\rho_A-\rho_{A'})}]
{}~~~\mbox{\footnotesize $A<A'$ for $L^+$, $A>A'$ for $L^-$}\\
{}\!\!\!\!\!{}\!\!\!\!\!&{}\!\!\!\!\!{}\!\!\!\!\! &
\Lpm{A}{B} (\T{A'}{B'}) = \mp \la \epsilon_A \epsilon_B + O(\la^2)
{}~~\mbox{\footnotesize $A\not=B,
A'\not= B$;~~$A<B$ for $L^+$, $A>B$ for $L^-$}\label{pmlanda2}
\ena
\noi all other $L^{\pm}(T)$ vanishing. Relations (\ref{pmlanda1})
and (\ref{pmlanda2}) imply that for any generator $\T{C}{D}$ we have
$
\Lpm{A}{B}(\T{C}{D}) = -\epsilon_A \epsilon_B \Lpm{B'}{A'}(\T{C}{D})
+
O(\la^2)~
\mbox{ with }\mbox{\footnotesize ~~ $A\not=B$, $A\not= B'$}~.$\\
In general, since
$$\Delta (\Lpm{A}{A})=\Lpm{A}{A} \otimes \Lpm{A}{A} ~~;~~~
\Delta(\Lpm{A}{B})=\Lpm{A}{A} \otimes \Lpm{A}{B}
+ \Lpm{A}{B} \otimes \Lpm{B}{B} + O(\la{}^2), ~~~\mbox{\footnotesize
$A\not=
B$}
$$
\noi we find that
\eqa
& &\Lpm{A}{A}=O(1)\\
& &\Lpm{A}{B}= O(\la),~~~~~\mbox{\footnotesize $A\not= B$,~$A\not=
B'$}\\
& &\Lpm{A}{A'}=O(\la^2) \mbox{~for $SO_q$},~~~~O(\la) \mbox{~for
$Sp_q$}
\ena
where, by definition, $\phi =O(\la{}^n)$ ($\phi$ being a functional)
means that
for any element $a\in \SqrNtwo$ with well-defined classical limit, we
have $\phi(a)=O(\la{}^n)$.

\noi Moreover the following relations hold:
\eqa
& &\Lpm{A}{A} = \Lmp{A}{A} + O(\la)~,\\
& &\kappa (\Lpm{A}{B})=\epsilon_A\epsilon_B\Lpm{B'}{A'}+O(\la) ~
\mbox{ and therefore,}~ \kappa^2=id + O(\la) ~.
\ena
Similarly one can prove the relations involving the
$f$ functionals (no sum on repeated indices):
\eqa
& &\ff{A}{AA}{A}=\epsi+O(\la)\label{fpropepsi}\\
& &\ff{A}{BA}{B}=O(1)~~~~ \mbox{ and
}~\ff{A}{BA}{B}=\ff{B'}{A'B'}{A'}+O(\la)
\label{relationf}\\
& &\ff{C}{CA}{A}=O(\la^2)~~~~\mbox{\footnotesize $C\not=
A$}\label{relationfCCAA}\\
& &\ff{C}{CA}{B}=O(\la^2)~~~~\mbox{\footnotesize [$A<B,C\not= B$] or
[$A>B,C\not= A$]}\label{6.30}
\ena
[hint: check (\ref{relationf})-(\ref{6.30}) first on the generators,
then use the coproduct in (\ref{copf})].
\noi From the last relation we deduce
\eqa
& &\cchi{A}{B}= \linv \ff{B}{BA}{B}+O(\la),~~~~~\mbox{\footnotesize
$A<B$}\\
& &\cchi{A}{B}= \linv \ff{A}{AA}{B}+O(\la),~~~~~\mbox{\footnotesize
$A>B$}
\ena
\noi and from (\ref{fpropepsi}) and (\ref{relationfCCAA}) one has
\eq
\cchi{A}{A}=\linv [\ff{A}{AA}{A}-\epsi]
\en
\noi Next one can verify that
\eq
\left. \begin{array}{ll} &\cchi{A}{B}(\T{B}{A})=-q_{BA}+O(\la)\\
                    & \cchi{A}{B}(\T{A'}{B'})=\epsilon_A
\epsilon_B+O(\la)\\
                    & \cchi{A}{B}(\T{C}{D})=0~ \mbox{ otherwise}
         \end{array}
                                   {}~~~~~~~~ \right\}
{}~\mbox{\footnotesize $A\not= B$, $A\not= B'$} \label{relationschiT}
\en
\eq
\forall  
~\T{C}{D}\;,~~~\cchi{A}{A}(\T{C}{D})=-\cchi{A'}{A'}(\T{C}{D})+O(\la).
\label{relationchiAAT}
\en
Eq.s (\ref{relationschiT}) yield the relation between $\chi$
functionals:

\eq
\forall ~\T{C}{D}\;,~~~\cchi{B'}{A'}(\T{C}{D})=-{\epsilon_A
\epsilon_B \over
q_{BA}}
\cchi{A}{B}(\T{C}{D})+O(\la),~
\sma{$A \not= B$, $A \not= B'.$}
\label{relationchiABT}
\en

It is not difficult to prove that the coproduct rule in (\ref{copf})
is compatible with
(\ref{relationchiABT}) and (\ref{relationchiAAT})
 making them valid on arbitrary polynomials in the $\T{A}{B}$
elements:
\eq
\cchi{B'}{A'}=-{\epsilon_A \epsilon_B \over q_{BA}}
\cchi{A}{B}+O(\la),~
\sma{$A \not= B$, $A \not= B'$}~~~;~~~~
\cchi{A}{A}=-\cchi{A'}{A'}+O(\la)~.\label{relationchiAA}
\en
Finally:
\eq
\cchi{A}{A'}=O(\la) \mbox{~for $SO_q$}~,~~~~O(1) \mbox{~for $Sp_q$}~,
 ~~~\sma{$A \not= A'$}. \label{relationchiAAp}
\en
Summarizing, in the $r \rightarrow 1$ limit, only the following
$\chi$
functionals survive:
\eqa
& &\cchi{A}{A}\equiv\limrone \linv~
[\ff{A}{AA}{A}-\epsi]\label{limuno}\\
& &\cchi{A}{B}\equiv\limrone \linv
{}~\ff{A}{AA}{B},~~~\sma{$A>B,A\not=B'$}\\
& &\cchi{A}{B}\equiv\limrone \linv ~
\ff{B}{BA}{B},~~~\sma{$A<B,A\not=B'$}\\
& &\cchi{A}{A'}\equiv\limrone \linv ~\sum_C \ff{C}{CA}{A'}=0
{}~~\sma{for
$SO_q$},
{}~~\not= 0~~\sma{for  $Sp_q$}\label{limquattro}
\ena
\noi Notice that  (\ref{relationchiAA}) and
 (\ref{relationchiAAp}) are all contained in the formula:
\eq
\cchi{B'}{A'}=-{\epsilon_A \epsilon_B \over q_{BA}}
\cchi{A}{B}+O(\la)\label{relationcchi}
\en
thus in the $\rone$ limit there are $\sma{$(N+2)(N+1)/2$}$ tangent
vectors for
$SO_q(N+2)$
and $\sma{$(N+2)(N+3)/2$}$ tangent vectors for $Sp_q(N+2)$, exactly
as in the
classical case.
\sk
The $r=1$ limit of  (\ref{dTAB}) reads:
\eq
d \T{A}{B}= -\sum_{C}\T{A}{C}  q_{CB} (\ome{B}{C}
- \epsilon_B\epsilon_C q_{BC}  \ome{C'}{B'})~, \label{dTABrone}
\en
and therefore, for $r=1$,  $\om$ appears only in the combination
\eq
\Ome{A}{B} \equiv \ome{A}{B} - \epsilon_A\epsilon_B q_{AB}
\ome{B'}{A'}~,
\label{defOme}
\en
Only $\sma{$(N+2)(N+1)/2$}$ [$\sma{$(N+2)(N+3)/2$}$ for $Sp_q(N+2)$]
of the $\sma{$(N+2)^2$}$ one forms $\Ome{A}{B}$
are linearly independent because [compare with (\ref{relationcchi})]:
\eq
\Ome{B'}{A'}=-{\epsilon_A\epsilon_B\over
q_{AB}}\Ome{A}{B}\label{relationOme}~.
\en
In the sequel, instead of considering the left module of one-forms
freely
generated by $\ome{A}{B}$, we consider the submodule $\Gamma$
freely generated by $\Ome{A}{B}$ with \sma{$A'<B$} for $SO_q$ and
 \sma{$A'\leq B$} for $Sp_q$.
In fact only this submodule will be relevant for the $r=1$
differential
calculus.
As in the classical case
 \footnote{
To make closer contact with the classical case one may define:
\[
\Omega^{AB} \equiv\Ome{C}{B} C^{CA}
=\epsilon_A\Ome{A'}{B}~~;~~~
\chil{AB} \equiv C_{AC} \cchi{C}{B} =\epsilon_A \cchi{A'}{B}  
~,\nonumber
\]
and  retrieve  the more familiar  $q$-antisymmetry:
\[
\Omega^{AB}=-\epsilon q_{BA} \Omega^{BA}~~;~~~
\chil{AB}= - \epsilon q_{AB} \chil{BA}~.\nonumber
\] } , in order to simplify notations in sums
we often use $\cchi{A}{B}$ and $\Ome{A}{B}$
without the restriction \sma{$A' \leq B$} see for ex. (\ref{defd3})  
below.
The bimodule structure on $\Gamma$, see {\sl Theorem 5.1}, is given  
by
the $\rone$ limit of the $\f{i}{j}$ functionals. These are diagonal  
in the
$i,j$
indices [i.e. they vanish for $i \not= j$, see
(\ref{fpropepsi})-(\ref{6.30})] and still satisfy the property  
(\ref{propf1}).
We have:
\eq
\begin{array}{rcl}
\Ome{A}{B}\,a &=&
(\ome{A}{B}-\epsilon_A\epsilon_B
q_{AB}\ome{B'}{A'})a\\
&=&(\ff{A}{BC}{D}*a)\ome{C}{D}-
\epsilon_A\epsilon_B q_{AB}(\ff{B'}{A'D'}{C'}*a)\ome{D'}{C'}\\
&=&(\ff{A}{BA}{B}*a)\Ome{A}{B}~
\label{bimOme}
\end{array}
\en
where in the last equality we have  used (\ref{relationf}) and no sum
is
understood.
We see that the bimodule structure is very simple since it does not
mix
different
$\Omega$'s. Moreover, relation (\ref{dTABrone}) is invertible and
yields:
\eq
\Ome{A}{B}=-q_{AB} \kappa (\T{B}{C}) d \T{C}{A} \label{Ome}~;
\en
in the limit $q_{AB}=1$, the $ \Ome{A}{B}$ are to be identified with
the
classical
one-forms, and indeed for $q_{AB}=1$ eq. (\ref{Ome}) reproduces
the correct classical limit $\Omega=-g^{-1}dg$ for the left-invariant
one-forms on the group manifold.
\sk
The bimodule commutation rule (\ref{bimOme}) yields a formula
similar to (\ref{commomT}). Replacing the values of the $R$ matrix
for $r=1$ we find the commutations:
\eq
 \Ome{A_1}{A_2} \T{R}{S}={q_{A_2 S} \over q_{A_1 S}} \T{R}{S}
 \Ome{A_1}{A_2} \label{commomTRone}
\en

For $r=1$ the coproduct on the $\chi$ functionals
reads
\eq
\Delta'(\cchi{A}{B})=\cchi{A}{B}\otimes\ff{A}{BA}{B} +
\epsi\otimes \cchi{A}{B} ~~~~~~~\mbox{ no sum on repeated
indices.}\label{copchiAB}
\en
cf. (\ref{copf}).
We then consider the $r=1$ limit of (\ref{defd2}) and therefore  
define  the
exterior differential
by:
\eq
da\equiv{1\over 2}(\cchi{A}{B}*a)\Ome{A}{B}=
\sum_{A'\leq B}(\cchi{A}{B}*a)\Ome{A}{B},~~~~~~\forall \/a\in A~,
\label{defd3}
\en
where in the second expression we have
used the basis of linear independent tangent
vectors  $\{\cchi{A}{B}\}_{A'\leq B}$
and dual one-forms $\{\Ome{A}{B}\}_{A'\leq B}$
(notice that in the $SO_q$ case we have
$\sma{$A'< B$}$ because $\cchi{A}{A}=\Ome{A}{A}=0$). The Leibniz
 rule is satisfied
for $d$ defined in (\ref{defd3}) because of (\ref{copchiAB}) and
(\ref{bimOme}).
Moreover any $\rho=a^A{}_B\Ome{A}{B}\in\Gamma$
can be written as $\rho=\sum_k a_kdb_k$, [use (\ref{Ome})].

We now introduce a left and a right action on the bimodule $\Gamma$
of
one-forms:
\eqa
&&\DL (a\Ome{A}{B})\equiv\D(a)I\otimes\Ome{A}{B}~,\label{defDL}\\
&&\DR (a\Ome{A}{B})\equiv\D(a)(\Ome{C}{D}\otimes
\MM{C}{DA}{B})~.\label{defDR3}
\ena
where $\MM{C}{DA}{B}=\T{C}{A}\kappa(\T{B}{D})$. [Using (\ref{defOme})
one can check that this is the $r=1$ limit of (\ref{deltaLomega}) and
(\ref{deltaRomega}).] Relation (\ref{defDR3}) is well defined i.e.
$ \DR(\Ome{B'}{A'})=\DR(-{\epsilon_A\epsilon_B\over
q_{AB}}\Ome{A}{B})$ because  
$\epsilon_F\epsilon_Eq_{FE}\MM{F'}{E'A}{B}=
\epsilon_A\epsilon_Bq_{AB}\MM{E}{FB'}{A'}$.
Since in the $r=1$ case the bicovariant bimodule conditions
(\ref{propf1}),
(\ref{copM}) and (\ref{propM}) are still satisfied, it is easy to  
deduce
 that $\DL$ and $\DR$ give
a bicovariant bimodule structure to $\Gamma$.
\sk
The differential (\ref{defd3}) gives a bicovariant differential
calculus if it
is
compatible with $\DL$ and $\DR$, i.e. if:
\eqa
&&\DL (adb)=\D(a)(id\otimes d)\D(b)~,\label{propdaRunoa}\\
&&\DR (adb)=\D(a)(d\otimes id)\D(b)~.\label{propdaRuno}
\ena
The proof of the compatibility of $d$ with $\DL$ is straightforward,
just use
(\ref{defd3}) and the coassociativity of the coproduct $\D$. In order
to prove
(\ref{propdaRuno}) it is sufficient to prove the following
\sk
\noi{\sl Theorem 6.1:}
\eq
\DR (db)=(d\otimes id)\D(b)~~~~~~~\forall ~ b ~.\label{propdaRsim}
\en
{\sl Proof: }$~$ We first review how the theorem is proved in the $r
\not= 1$
case.
On the left hand side we have
\eq
\DR (db)=\DR[(\chi_i*b)\om^i]=\D[b_1\chi_i(b_2)]\om^j\otimes\M{j}{i}
=b_1\om^j\otimes b_2\chi_i(b_3)\M{j}{i} \label{p1}
\en
\noi with $(\D \otimes id)\D (b)=
(id \otimes \D) \D (b) \equiv b_1 \otimes b_2 \otimes b_3$
 [cf. (\ref{axiom1})],  while for the right hand side
\eq
(d\otimes id)\D (b)=db_1 \otimes b_2=
b_1\chi_j(b_2)\om^j
\otimes b_3=b_1\om^j\otimes\chi_j(b_2)b_3~. \label{p2}
\en
Therefore (\ref{propdaRsim}) holds if and only if
\eq
b_1\om^j\otimes b_2 \chi_i
(b_3)\M{j}{i}=b_1\om^j\otimes\chi_j(b_2)b_3\label{OmechiM}
\en
and this last relation is equivalent to
\eq
b_1\chi_i( b_2)\M{j}{i}=\chi_j(b_1)b_2~,
{}~~\mbox{ i.e. }~~~\chi_i*b=(b*\chi_j)\kappa(\M{i}{j})   
\label{leftrightchi}
\en
as one can verify by applying $m(\kappa \otimes id) \DL \otimes id $
($m$ denotes multiplication) to (\ref{OmechiM}), and using the linear
independence of the $\om^i$. Now  formula (\ref{leftrightchi}) holds
also in the limit $r=1$.  Indeed if we
consider $b$ to be a polynomial in the $\T{A}{B}$ with well behaved
coefficents
in the
$\rone$ limit, then
$\lim_{\rone}[b_1\chi_i(b_2)\M{j}{i}]=\lim_{\rone}[\chi_j(b_1)b_2]$
i.e.
$b_1[\lim_{\rone}\chi_i(b_2)]\M{j}{i}=
[\lim_{\rone}\chi_j(b_1)]b_2$ so that relation (\ref{leftrightchi})
remains valid  for $r=1$, cf .(\ref{limuno})-(\ref{limquattro}).
At this point one can prove {\sl Theorem 6.1} in the $r=1$ case
simply by substituting $\Omega$ to $\om$ in (\ref{p1}), (\ref{p2})
and
(\ref{OmechiM}). Since (\ref{leftrightchi}) holds for $r=1$, then
also (\ref{OmechiM}) holds in this limit and the theorem is proved.
{}~~~~$\Box$
\sk
Relation (\ref{leftrightchi})  has an
important geometrical interpretation: the left invariant vector field
$\chi_i*$
(associated with the tangent vector $\chi_i$)  can be expressed in
terms
of the right invariant vector fields $*\chi_i$ via the
``deformed functions on the group'' $\kappa(\M{j}{i})$
\cite{PaoloPeter}.
\sk
In virtue of {\sl Remark 5.1} we conclude that (\ref{defd3}) defines
a bicovariant differential calculus on $S_q(N+2)$.
\sk
\noi{\sl Note 6.1 :}
in the right-hand side of  (\ref{leftrightchi}) the sum on
the indices $j=(\mbox{\sma{$C,D$}})$ can be restricted to \sma{$C'  
\leq D$} ,
thus using the basis $\{\cchi{C}{D}\}_{C'\leq D}$,
provided
one replaces $M$ by
\eq
\begin{array}{lll}
\MMc{C}{DA}{B} \equiv \MM{C}{DA}{B}-\epsilon_C \epsilon_D
q_{DC}
\MM{D'}{C'A}{B} ~& &~\sma{$\mbox{ for } ~C'
\not= D ,~A'\not= B$} \\
\MMc{C}{C' A}{B} \equiv 0 ~&
\!\!\!\!\!\!\!\!\!\!\!\!\!\!\!\!\!\!\!\!\!\!\!\!\!\!\!\!\!\!\!\!\!\!\! 
\!\!\!\!\!\!\!\!\!\!\!\!\!\!\!\!\!\!\!
,~~\MMc{C}{D A}{A'} \equiv 0 &~\sma{$\mbox{ for } SO_q$} \\
\MMc{C}{C' A}{B}\equiv\MM{C}{C' A}{B} ~&
\!\!\!\!\!\!\!\!\!\!\!\!\!\!\!\!\!\!\!\!\!\!\!\!\!\!\!\!\!\!\!\!\!\!\! 
\!\!\!\!\!\!\!\!\!\!\!\!\!\!\!\!\!\!\!
,~~\MMc{C}{D A}{A'}\equiv\MM{C}{D A}{A'}
 &~\sma{$\mbox{ for } Sp_q$}\label{Mantisymm}
\end{array}
\en
This is easily seen from
(\ref{relationcchi}).
We can also write
$\DR (a\Ome{A}{B})=\sum_{C'\leq D}\D(a)(\Ome{C}{D}\otimes
\MMc{C}{DA}{B})~$  cf.(\ref{defDR3}), thus using the basis
$\{\Ome{C}{D}\}_{C'\leq D}$. According to the general theory  
\cite{Wor}, the
elements
$\MMc{C}{DA}{B}$ with \sma{$C' \leq D,~A' \leq B$} are by definition  
the
adjoint elements for the differential calculus on
$S_{q,r=1}(N+2)$.
Since the calculus is bicovariant [cf.(\ref{propdaRunoa}),  
(\ref{propdaRuno})]
we  know a priori that
the $\MMc{C}{DA}{B}$ with \sma{$C' \leq D,~A' \leq B$} satisfy the
properties
(\ref{copM}) and (\ref{propM}).  \footnote{
A direct proof in the \sma{$ SO_q$} case is also instructive.
We call $P_{\!-}$ the ``q-antisymmetric" projector
defined by:
\[
P_{\!-}{}^{A~~~D}_{~BC}\equiv\onehalf(\de^A_C\de^D_B-q_{BA}\de_C^{B'}\ 
de_{A'}^D)=\onehalf(\de^A_C\de^D_B-q_{CD}\de_C^{B'}\de_{A'}^D)~.
\]
Then one easily shows that
$~P_{\!-}{}^{A~~~D}_{~BC}=-q_{BA}P_{\!-}{}^{B'~~~C}_{~A'D}~
,~~ P_{\!-}{}^{A~~~D}_{~BC}=-q_{CD}P_{\!-}{}^{A~~~C'}_{~BD'}$ and
\[
\begin{array}{c}
\Omega^jP_{\!-}{}_j{}^i=\Omega^i ~,~~~P_{\!-}{}_i{}^j\chi_j=\chi_i  
~,~~~
P_{\!-}{}_k{}^if^k{}_j=f^i{}_kP_{\!-}{}_j{}^k=P_{\!-}{}_k{}^if^k{}_nP_ 
{\!-}{}_j{}^n~,\\
\MMcc{i}{j}=2P_{\!-}{}_i{}^lM_{l}{}^{j}=2M_{i}{}^{l}P_{\!-}{}_l{}^j~,~ 
~~
\MMcc{i}{j}=P_{\!-}{}_i{}^l\MMcc{l}{j}=\MMcc{i}{l}P_{\!-}{}_l{}^j=
2P_{\!-}{}_i{}^{\alpha}\MMcc{\alpha}{j}
=2\MMcc{i}{\beta}P_{\!-}{}_{\beta}{}^{j}~,\\
\end{array}
\]
where greek letters $\alpha, \beta$ represent adjoint indices
\sma{$(A_1,A_2),~(B_1,B_2)$} with the restriction
\sma{$A_1' < A_2,~B_1' < B_2$}. It is then straightforward to show  
that
$\D(\MMcc{i}{j})=\MMcc{i}{\alpha}\otimes\MMcc{\alpha}{j}$ and
$\epsi(\MMcc{\alpha}{\beta})=\de_{\alpha}^{\beta}.$ Applying  
$P_{\!-}$ to
(\ref{propM}) and
using $f^i{}_j=0$ unless $i=j$ cf.(\ref{fpropepsi})-(\ref{6.30}) one  
also
proves
$\M{\alpha}{j} (a * \f{\alpha}{k})=(\f{j}{\beta} * a) \M{k}{\beta}$.  
These
formulae
hold in particular if all indices are greek, thus proving
(\ref{copM}) and (\ref{propM}) for $SO_{q,r=1}(N+2)$.}

\sk
It is useful to express the bicovariant
algebra
(\ref{bico1}), (\ref{bico2})-(\ref{bico4})  in the $\rone$ limit.
Due to the $R$ matrix being diagonal for $r=1$, the $\Lambda$ tensor
$
\LL{A_1}{A_2}{D_1}{D_2}{C_1}{C_2}{B_1}{B_2}
\equiv \ff{A_1}{A_2B_1}{B_2} (\MM{C_1}{C_2D_1}{D_2})
$
takes the simple form:
\eq
\LL{A_1}{A_2}{B_1}{B_2}{B_1}{B_2}{A_1}{A_2}=q_{A_1B_2}q_{A_2B_1}
q_{B_1A_1}q_{B_2A_2}~,
{}~~~~0 \mbox{ otherwise}
\en
Therefore (\ref{bico2})-(\ref{bico4}) read (no sum on repeated
indices):
\eq
\f{i}{i}\f{j}{j}=\f{j}{j}\f{i}{i}
\en
\eq
\C{jk}{i} \f{j}{j} \f{k}{k} + \f{i}{j} \chi_k= \L{kj}{jk} \chi_k
\f{i}{j} + \C{jk}{i} \f{i}{i} \label{1bico3}
\en
\eq
\chi_k  \f{i}{i}=\L{ik}{ki} \f{i}{i} \chi_k~.  \label{1bico4}
\en
Explicitly the $q$-Lie algebra (\ref{bico1}) reads:
\eqa
& &\cchi{C_{1}}{C_2}  \cchi{B_{1}}{B_2} - q_{B_1C_2} q_{C_1B_1}
q_{B_2C_1} q_{C_2B_2}
{}~\cchi{B_{1}}{B_2} \cchi{C_{1}}{C_2}=~~~~~~~~~\nonumber\\
& &~~~~-q_{B_1C_2} q_{C_2B_2} q_{B_2B_1} \de^{C_1}_{B_2}
{}~\cchi{B_{1}}
{C_2}+
q_{C_1B_1} q_{B_2B_1} C_{B_2C_2} ~\cchi{B_{1}}{C_1'}+\nonumber\\
& &~~~~+q_{C_2B_2} q_{B_1C_2} C^{C_1B_1}~ \cchi{B_2'} {C_2}-
q_{B_2C_1} \de^{B_1}_{C_2}~\cchi{B_2'}{C_1'}~. \label{Lierone}
\ena
The Cartan-Maurer equations are obtained by differentiating
(\ref{Ome}):
\eq
d\Ome{A}{B}=q_{AB} q_{BC} q_{CA} C_{CD} ~\Ome{C}{B} \we \Ome{A}{D}
\label{CMrone}
\en

The commutations between $\Omega$ 's are easy to find using
(\ref{omcomrone}):
\eq
\Ome{A_1}{A_2} \we \Ome{D_1}{D_2} = - q_{A_1 D_2} q_{D_1A_1}
q_{A_2D_1} q_{D_2A_2} \Ome{D_1}{D_2} \we \Ome{A_1}{A_2}
\label{commOm}
\en
\sk

Finally, we turn to the ${}^*$-conjugations given by equations
(\ref{conjchiAB}) and
(\ref{conjomAB}).
Their $\rone$ limit yields
\eq
(\Ome{A}{B})^*=-q_{BA} \Dcal{C}{A} \Ome{C}{D} \Dcal{B}{D}~~~;~~~~
(\cchi{A}{B})^*=-q_{CD} \Dcal{A}{C} \cchi{C}{D} \Dcal{D}{B}~,
\label{conjchiABo}
\en
and shows that we have a bicovariant ${}^*$-differential calculus.


\sect{Differential calculus on $\ISOqroN$}


We reconsider now, in the $\rone$ limit,  the functionals given
in eq.s (\ref{chiN2})-(\ref{chiN2end}). We list below the functionals
among
these that
annihilate the Hopf ideal:
\eqa
& & \cchi{a}{b}=\rinv [\ff{c}{c a}{b}-\de^a_b \epsi] \nonumber \\
& & \cchi{a}{\ci}=\rinv \ff{c}{c a}{\ci} \nonumber \\
& & \cchi{\bu}{b}=\rinv \ff{\bu}{\bu\bu}{b} \nonumber \\
& &~~\nonumber \\
& &\cchi{\ci}{\ci}=\rinv [\ff{\ci}{\ci\ci}{\ci}-\epsi]  \nonumber \\
& & \cchi{\bu}{\bu} = \rinv [\ff{\bu}{\bu\bu}{\bu} - \epsi]
\label{chiin}
\ena
\noi Note that in the $\rone$ limit
 $\cchi{\bu}{\ci}$ vanishes  for
 $\SOqroNt$, and does not vanish in the case $\SpqroNt$.
We treat here the orthogonal case, both for simplicity and because we
are more interested for physical reasons
to orthogonal (rather than symplectic) inhomogeneous q-groups.
The reader can easily extend our discussion to the symplectic case,
and include the $\cchi{\bu}{\ci}$ tangent vector (besides the
diagonal $\cchi{a'}{a}$).
\sk
Taking all the $\chi$'s given in (\ref{chiin}) one obtains
 a differential calculus
containing dilatations (because of the presence of $\cchi{\ci}{\ci}$
and $\cchi{\bu}{\bu}$). It is however possible to exclude
the generators  $\cchi{\ci}{\ci}$
and $\cchi{\bu}{\bu}$ from the list,  and obtain a
dilatation-free
bicovariant differential calculus. This we will discuss in the rest
of this section,
while, in a more general setting, the  case with dilatations is  
discussed in
ref.
\cite{UEA}.
\sk
For $r=1$ the $\chi$'s in (\ref{chiin}) are not independent, cf.
relation (\ref{relationcchi}) of previous section, and we have:
\eq
\cchi{b'}{a'}=-q_{ab} \cchi{a}{b},~~\cchi{b'}{\ci}=
-{1\over q_{b \bu}} \cchi{\bu}{b},~~\cchi{\ci}{\ci}=-\cchi{\bu}{\bu}
\en
Therefore we consider $\cchi{a}{b}\/(a'<b), \cchi{\bu}{b}$ as a
candidate basis for the tangent vectors on
$\ISOqroN$. We will show  that these
$\chi$ functionals indeed define a bicovariant differential calculus
on $\ISOqroN$.
\sk
{\sl Theorem 7.1:} the functionals $\f{i}{j}$, obtained from those of
$\SOqroNt$ by restricting  the indices
to $i=ab,\bu b$, annihilate the Hopf ideal $H$.

{\sl Proof:}
According to the results of the previous section,
the only non-vanishing functionals with indices
$i=ab,\bu b$ are
\eqa
& &\ff{a_1}{a_2b_1}{b_2}= \kappa (\Lp{b_1}{a_1}) \Lm{a_2}{b_2}
\nonumber\\
& &\ff{\bu}{a_2 \bu}{b_2}= \kappa (\Lp{\bu}{\bu}) \Lm{a_2}{b_2}
\label{fin}
\ena
\noi To prove the theorem, first one checks directly that the
functionals (\ref{fin}) vanish
on the generators ${\cal T}$ of the ideal  $H$, i.e. on ${\cal T} =
\T{a}{\circ}, \T{\bu}{b}, \T{\bu}{\ci}$.
This extends to any element
of the form $a {\cal T} b$ ($a, b  \in  \SOqroNt$), i.e. to any
element of
$H=Ker(P)$, because of the property (\ref{propf1}) which in the
$\SOqroNt$ reads $\Dtwo '(f^i{}_i)=f^i{}_i\otimes f^i{}_i$ since
the functionals $\f{i}{j}$ vanish when $i\not= j$.  \square
\sk
Thus the functionals $\chi_i$ and $\f{i}{i}$ with  $i=ab,\bu  b$,
which we denote collectively by the symbol $f$,
all vanish on $H$.  Then these functionals are
well defined on the quotient $\ISOqroN = \SOqroNt / Ker(H)$,
in the sense that the ``projected" functionals
\eq
\fb : \ISOqroN \rightarrow {\Cb},~~
\fb (P(a)) \equiv f(a)~ , ~~~\forall a \in \SOqroNt \label{deffb}
\en
are well defined.
Indeed if $P(a)=P(b)$, then $f(a)=f(b)$ because
$f(Ker(P))=0$. This holds for any
 functional $f$ vanishing on $Ker(P)$.
\sk
The product $fg$ of  two generic functionals
vanishing on $KerP$ also vanishes
on $KerP$, because $KerP$ is a co-ideal (see ref. \cite{qpoincare1}):
 $fg(KerP)=(f\otimes g)\D_{N+2} (KerP)=0$.
Therefore $\overline{fg}$ is well defined on $\ISOqroN$; moreover,  
[use
(\ref{co-iso})]
$\overline{fg}[P(a)]\equiv fg(a)= (\bar{f}\otimes \bar{g}) \D(P(a))  
\equiv
 \bar{f} \bar{g}[P(a)] $,
\noi and the product of $\bar{f}$ and $\bar{g}$ involves
the coproduct $\D$ of $\ISOqroN$.
\sk
There is a natural way to introduce a coproduct on the $\fb$'s :
\eq
\D '\fb[P(a)\otimes P(b)]\equiv\fb[P(a)P(b)]=\fb[P(ab)]
=f(ab)=\D_{N+2}'(f)[a\otimes b] ~.
\en
It is then straightforward to show, from the relations (\ref{copf})  
for
$\SOqroNt$ i.e.
$\Dtwo '\chi_i=\chi_i\otimes f^i{}_i +\epsi\otimes\chi_i ,~\Dtwo
'f^i{}_i=f^i{}_i\otimes f^i{}_i$, that
\eqa
& &\D '\fb^{i}{}_{i} = \fb^{i}{}_{i}\otimes \fb^{i}{}_{i} ~~\mbox{  
i.e. }~
\fb^{i}{}_{i}[P(a)P(b)] = \fb^{i}{}_{i}[P(a)]\fb^{i}{}_{i}[P(b)]  
\label{cofa}\\
& &\D '{\bar \chi}_{i}={\bar \chi}_{i} \otimes \fb^{i}{}_{i} +
 \epsi \otimes {\bar \chi}_{i}
\label{cofb}
\ena
with $i$ adjoint index running over the set of indices
$ab,\bu b$.
With abuse of notations we
will simply write $f$ instead of $\fb$,  and the
$f$  in (\ref{fin}) will be seen as functionals on
$\ISOqroN$.
\sk
Consider now the elements $\Mc{i}{j} \in \ISOqroN$ obtained by  
projecting with
$P$
those of  $\SOqroNt$ and
with the restriction $i=ab,\bu b$ on the adjoint
indices.
The effect of the projection is to replace
the coinverse in $\SOqroNt$, i.e. $\kappa_{N+2}~,$
with the coinverse $\kappa$ of $\ISqroN$ (see the
last of (\ref{co-iso})). The nonvanishing elements are:
\eqa
& &\MMc{b_1}{b_2a_1}{a_2} = \T{b_1}{a_1} \kappa (\T{a_2}{b_2})
- q_{b_2b_1} \T{b_2'}{a_1} \kappa (\T{a_2}{b_1'})\nonumber\\
& &\MMc{b_1}{b_2 \bu}{a_2} = x^{b_1} \kappa (\T{a_2}{b_2})
- q_{b_2b_1} x^{b_2'} \kappa (\T{a_2}{b_1'})
\nonumber\\
& &\MMc{\bu}{b_2 \bu}{a_2} = v  \kappa (\T{a_2}{b_2})
\label{Min}
\ena

In the sequel  greek  letters will denote adjoint indices
$\al= (a_1,a_2)$ with $a'_1 < a_2 \mbox{, and } \al=(\bu,a_2)$.
\sk
{\sl Theorem 7.2:} the left $A$-module [$A=\ISOqroN$]
 $\Ga$ freely generated
by the symbols $\Omega^{\al}$
is a bicovariant bimodule over $\ISOqroN$
with the right multiplication (no sum on repeated indices):
\eq
\Omega^{\al} a = (\f{\al}{\al} * a) \Omega^{\al}~,
{}~~~~a \in \ISOqroN \label{omia}
\en
\noi [where the $\f{\al}{{\be}}$ are found in (\ref{fin}) and the  
$*$-product
is computed with the co-product $\D$ of $\ISOqroN$] and with
the left and right actions
of $\ISOqroN$ on $\Ga$  given by:
\eqa
& &\DL (a_{\al} \Omega^{\al}) \equiv \D(a_{\al}) I \otimes  
\Omega^{\al}
 \label{DLin}\\
& &\DR (a_{\al} \Omega^{\al}) \equiv \D(a_{\al}) \Omega^{\be} \otimes
\Mc{{\be}}{{\al}}
\label{DRin}
\ena
\noi the $\Mc{{\be}}{\al}$ being given in (\ref{Min}).
\sk
{\sl Proof:} we prove the theorem
by showing  that  the
functionals $f$
and the elements ${\Mc{\al}{\be}}$ listed in (\ref{fin}) and
(\ref{Min})
satisfy the properties (\ref{propf1})-(\ref{propM})
(cf. the theorem by Woronowicz discussed in the Section 4).
Applying the projection $P$ to the $\SOqroNt$ relation
$\Dtwo(\Mc{\al}{\be})=\Mc{\al}{\ga}\otimes\Mc{\ga}{\be}$,
one verifies directly that
the elements $\Mc{\al}{\be}$
in (\ref{Min}) do satisfy the properties
(\ref{copM}).
We have already shown that
the functionals $f$ in (\ref{fin}) satisfy
(\ref{propf1}).

Consider now the last property (\ref{propM}).
For $\SOqroNt$ it explicitly reads (cf.  {\sl Note 6.1} ):
\eq
\MMc{A_1}{A_2B_1}{B_2} (a *  
\ff{A_1}{A_2C_1}{C_2})=(\ff{B_1}{B_2A_1}{A_2} * a)
\MMc{C_1}{C_2A_1}{A_2}~~~
\mbox{ with \sma{$ A_1'<A_2,~B_1'<B_2,~C_1'<C_2 $}}\; .  
\label{propM7}
\en
Restrict the free indices
to greek indices, and
apply the projection
$P$ on
both members of the equation.
It is then  immediate
to see that  only the $f$'s in (\ref{fin}) and the
$M_{\!-}$'s in (\ref{Min}) enter in (\ref{propM7}).

We still have to prove that the $*$ product
in (\ref{propM7})
can be computed via the coproduct $\D$ in
$\ISOqroN$.
Consider the projection of property
 (\ref{propM7}), written symbolically as:
\eq
P [ {M_{\!-}} (f \otimes id) \D_{N+2} (a)]=P [(id\otimes f)
\D_{N+2} (a) M_{\!-}]~. \label{PMfa}
\en
Now apply the definition (\ref{deffb}) and the first of
(\ref{co-iso}) to rewrite
(\ref{PMfa}) as
\eq
P({M_{\!-}})(\fb\otimes  
id)\D(P(a))=(id\otimes\fb)\D(P(a))P({M_{\!-}})~.
\en
This projected equation then  {\sl becomes}
 property (\ref{propM})  for the $\ISOqroN$
functionals
$f$ and adjoint elements ${M_{\!-}}$, with the
correct coproduct $\D$ of $\ISOqroN$.  \square
\sk
Using the general formula (\ref{omia}) we can
deduce the $\Omega, T$ commutations for $\ISOqroN$:
\eqa
& &\Ome{a_1}{a_2} \T{r}{s}={q_{a_2 s} \over q_{a_1 s}} \T{r}{s}
 \Ome{a_1}{a_2}\\
& &\Ome{a_1}{a_2} x^r= {q_{a_2 \bu} \over q_{a_1 \bu}} x^r
\Ome{a_1}{a_2} \\
& &\Ome{a_1}{a_2} u= {q_{a_1 \bu} \over q_{a_2 \bu}} u~
\Ome{a_1}{a_2} \\
& &\Ome{\bu}{a_2} \T{r}{s}=q_{s\bu} q_{a_2 s} \T{r}{s}
\Ome{\bu}{a_2}\label{VTcomm}\\
& &\Ome{\bu}{a_2} x^r=q_{a_2\bu} x^r \Ome{\bu}{a_2} \label{Vxcomm}\\
& &\Ome{\bu}{a_2} u={1\over q_{a_2\bu}} u~\Ome{\bu}{a_2}
\ena
\noi {\sl Note 7.1:} $u$ commutes with all
 $\Omega$ 's only if $q_{a\bu}=1$ (cf. {\sl Note 3.2}). This means  
that
$u=I$ is consistent with the differential calculus
on $ISO_{q_{ab},r=1,q_{a\bu}=1}(N)$.
\sk
An exterior derivative on $\ISOqroN$ can be defined
 as
\eq
da=(\chi_{\al} * a) \Omega^{\al}
\en
where the $\chi_{\al} = \cchi{a}{b} \/ (a'<b), \cchi{\bu}{b}$ are  
given in
(\ref{chiin}). Due to the coproduct  (\ref{cofb}) and
the commutations (\ref{omia}) this derivative
satisfies the Leibniz rule. It is also compatible
with the left and right action of $ISO_{q,r=1}(N)$
since (\ref{propdaRuno}) holds. This can be
seen by noting that the key equation (\ref{leftrightchi}), which
in the $\SOqroNt$ case reads [see  
(\ref{Mantisymm}),~(\ref{relationcchi})]:
\eq
\cchi{A_1}{A_2}*b=(b*\cchi{B_1}{B_2})
\kappa(\MMc{A_1}{A_2B_1}{B_2})~~~
\mbox{ with \sma{$ A_1'<A_2,~B_1'<B_2$} }
\en
becomes
property (\ref{leftrightchi})  for the $\ISOqroN$
functionals
$\chi$ and adjoint elements ${M_{\!-}}$, with the
correct coproduct $\D$ of $\ISOqroN$, once we restrict the free  
indices
to greek indices and apply the projection $P$.
\sk
The exterior
derivative on the generators $\T{A}{B}$ is given by:
\eqa
& & d \T{a}{b}= -\sum_{c}\T{a}{c} q_{cb} \Ome{b}{c}\\
& & dx^a=-\sum_c \T{a}{c} q_{c\bu} V^c\\
& & du=dv=0 \label{dTABiso}
\ena
\noi  where we have defined $V^a \equiv \Ome{\bu}{a}$.
Again, for $q_{a\bu}=1$,  $u=v=I$ is a consistent  choice.
\sk
Every element $\rho$ of $\Ga$ can be written as
$\rho=\sum_k a_k db_k$ for some $a_k,b_k$ belonging to
$\ISOqroN$. Indeed inverting (\ref{dTABiso}) yields:
\eqa
& & \Ome{a}{b}=-q_{ab} \kappa (\T{b}{c}) ~d \T{c}{a}\label{Omein}\\
& & V^b=-{1\over q_{b\bu}}
\kappa (\T{b}{c}) ~ dx^c \label{Vin}
\ena
\sk
Thus all the axioms for a bicovariant
first order differential calculus on $\ISOqroN$
are satisfied.
\sk
The exterior product of the left-invariant one-forms
is defined as
\eq
\Om^{\al} \we \Om^{\be}\equiv \Om^{\al} \otimes \Om^{\be} -  
\L{\al\be}{\ga\de}
\Om^{\ga} \otimes \Om^{\de}
\en
\noi where
\eq
\L{\al\be}{\ga\de}=\f{\al}{\de} (\Mc{\ga}{\be})
\en
This $\Lambda$ tensor can in fact
be obtained from
the one of $\SOqroNt$ by restricting its indices to the
subset $ab,  \bu b$. This is true because
when $i,l=ab, \bu b$  we have
$\f{i}{l} (Ker P)=0$ so that $\f{i}{l}$ is
well defined on $\ISOqroN$, and we can write
$\f{i}{l} (\Mc{k}{j})=\fb^i_{~l} [P(\Mc{k}{j})]$
(see discussion after
{\sl Theorem 7.1}). Then we can just specialize indices
in equation (\ref{commOm}) and deduce
the $q$-commutations for the one-forms $\Omega$ and $V$:
\eq
\Ome{a_1}{a_2} \we \Ome{d_1}{d_2} = - q_{a_1 d_2} q_{d_1a_1}
q_{a_2d_1} q_{d_2a_2} \Ome{d_1}{d_2} \we \Ome{a_1}{a_2}
\en
\eq
\Ome{a_1}{a_2} \we V^{d_2} = - {q_{a_2\bu} \over q_{a_1\bu}}
q_{a_1 d_2} q_{d_2a_2} V^{d_2} \we \Ome{a_1}{a_2}
\en
\eq
V^{a_2} \we V^{d_2}=-{q_{a_2\bu} \over q_{d_2\bu}} q_{d_2a_2}
V^{d_2} \we V^{a_2}
\en
\sk
The exterior differential on $\Ga^{\we n}$
can be defined as in Section 5 (eq. (\ref{defdonga})),
and satisfies all the properties (\ref{propd2})-(\ref{propd4}).
\sk
The Cartan-Maurer equations
\eq
d\Om^{\al}=
-\onehalf\c{\be\ga}{\al} \Om^{\be} \we \Om^{\ga}\label{CMin0}
\en
\noi can be explicitly written for the $\Omega$ and $V$
by differentiating eq.s (\ref{Omein}) and (\ref{Vin}):
\eqa
& &d\Ome{a}{b}= q_{ab} q_{bc} q_{ca} ~\Ome{c}{b} \we \Ome{a}{c}
\label{CMin1}\\
& &d V^b={q_{a\bu} \over q_{b\bu}} q_{ba} ~\Ome{a}{b} \we V^a
\label{CMin2}
\ena
where the one-forms
$\Ome{a}{b}$ with $a'>b$ are given by  $\Ome{a}{b}=-q_{ab}  
\Ome{b'}{a'}$;
i.e. we consider (as it is usually done in the classical limit),
the one-forms $\Ome{a}{b}$ to be
``$q$-antisymmetric''
$\Ome{a}{b}=-q_{ab} \Ome{b'}{a'}$,
cf. eq. (\ref{relationOme}).
\sk
Using the values of  $\L{\al\be}{\ga\de}=\f{\al}{\de}  
(\Mc{\ga}{\be})$
and of the structure constants $\C{\al\be}{\ga}=
\chi_k (\Mc{\be}{\al})$
we can explicitly write the ``$q$-Lie algebra"
of $\ISOqroN$. The $\cchi{c_1}{c_2}$, $\cchi{b_1}{b_2}$
$q$-commutations
read as in eq. (\ref{Lierone}) with lower case indices, and
give the $SO_{q,r=1}(N)$ $q$-Lie algebra; the remaining
commutations are
\eqa
& &\cchi{c_1}{c_2} \chi_{b_2} - {q_{c_1\bu} \over q_{c_2\bu}}
q_{b_2c_1} q_{c_2b_2} \chi_{b_2} \cchi{c_1}{c_2}=
{q_{c_1\bu} \over q_{c_2\bu}}[C_{b_2c_2} \chi_{c_1'}-\de^{c_1}_{b_2}
q_{c_2c_1} \chi_{c_2}] \nonumber\\
& &\chi_{c_2} \chi_{b_2} - {q_{b_2\bu} \over q_{c_2\bu}}
q_{c_2b_2} \chi_{b_2} \chi_{c_2} \label{qLie3}
=0
\ena
\noi with the definition
\eq
\chi_a \equiv  \cchi{\bu}{a}
\en
It is not difficult to verify that the $\Cb$ constants
 do coincide with the $C$
constants appearing in the Cartan-Maurer equations
(\ref{CMin0})-(\ref{CMin2}).
\sk

The *-conjugation on the $\chi$ functionals and on the one-forms
$\Omega$ can be deduced from (\ref{conjchiABo}):
\eq
(\cchi{a}{b})^*=-q_{cd} \Dcal{a}{c} \cchi{c}{d}
\Dcal{d}{b},~~(\chi_b)^*=
-(q_{d\bu})^{-1} \chi_d \Dcal{d}{b}
\en
\eq
(\Ome{a}{b})^*=-q_{ba} \Dcal{c}{a} \Ome{c}{d} \Dcal{b}{d},~~(V^b)^*=
-q_{b\bu} V^d \Dcal{b}{d}
\en
{\sl Remark 7.1:}  as discussed in \cite{qpoincare1}  and at end of
Section 3,
a  $q$-Poincar\'e group without dilatations (i.e. $u=I$) has only
one free real parameter $q_{12}$, which is the real para-
meter
related to the $q$-Lorentz subalgebra. Then the formulas of this
section
can be specialized to describe a bicovariant calculus on the
dilatation-free
$ISO_{q,r=1}(3,1)$ provided  $q_{a\bu}=1$ and  $q_{12} \in \Rb$. It
is
however possible to have a bicovariant calculus without the  
dilatation
generator
$\cchi{\bu}{\bu}$ even on $ISO_{q,r=1}(3,1)$ with $u \not= I$. The
possibility
of having a dilatation-free $q$-Lie algebra describing a bicovariant
calculus
on
a $q$-group containing dilatations $u$ was already observed in the
case of $IGL$ $q$-groups, see ref. \cite{qigl2}. The $q$-Poincar\'e
algebra presented in \cite{qpoincarediff1} corresponds to the case
$q \equiv q_{1\bu}$, $q_{2\bu}=q_{12}=1$, for which the Lorentz
subalgebra is undeformed and the $q$-Poincar\'e group contains  
$u\not= I$.
Finally, the bicovariant calculus that
includes the
dilatation generator $\cchi{\bu}{\bu}$ is discussed in ref.
\cite{UEA}.


\sect{The multiparametric orthogonal quantum plane as a $q$-coset
space}


In this section we derive the differential calculus
on the orthogonal quantum plane
\eq
Fun_{q,r=1}\left( {ISO(N) \over SO(N)}\right)~,
\en
\noi i.e. the $ISO_{q,r=1}(N)$ subalgebra
generated by the coordinates $x^a$.

The coordinates $x^a$
satisfy the commutations
(\ref{PRTT13}):
\eq
x^a x^b = q_{ab} ~x^b x^a
\en
\noi Note that the coordinates $x^a$ do not trivially
commute with the $SO_{q,r=1}(N)$ $q$-group elements, but
$q$-commute according to relations (\ref{PRTT33}):
\eq
\T{b}{d} x^a={q_{ba} \over q_{d\bullet}}  x^a \T{b}{d}
\label{PRTT33bis}
\en
\noi {\sl Lemma:} $\cchi{b}{c} (a)=0$ when $a$ is a
polynomial in $x^a$ and $v$, with all monomials
containing at least one $x^a$. This is easily proved
by observing that no tensor exists with the correct
index structure.  In fact we can extend this lemma
even to  $v \cdots v$, due to
\eq
\cchi{b}{c} (v)=0
\en
\noi and the coproduct rule (\ref{cofb}), reading explicitly (no sum  
on
repeated indices):
\eq
\begin{array}{lll}
& &\D '(\cchi{a}{b})=\cchi{a}{b} \otimes \ff{a}{ba}{b}+
\epsi \otimes \cchi{a}{b} ~,\\
& &\D '(\cchi{\bu}{b})=
\cchi{\bu}{b} \otimes \ff{\bu}{b\bu}{b} +
\epsi \otimes \cchi{\bu}{b}~.
\end{array} \label{copex}
\en
\sk
\noi {\sl Theorem 8.1:} $\cchi{b}{c} * a=0$ when $a$ is
a polynomial in $x^a$.
\sk
\noi {\sl Proof:} we have $\cchi{b}{c} * a=(id\otimes
\cchi{b}{c})(a_1\otimes
a_2)=a_1 \cchi{b}{c} (a_2)$. We use here the standard notation
$\D (a) \equiv a_1 \otimes a_2$.  Since $a_2$ is a polynomial
in $x^a$ and $v$ (use the coproduct rule (\ref{Pcoproduct2})), and
$\cchi{b}{c}$ vanishes on such a polynomial
(previous Lemma),
the theorem is proved. ~~~\square
\sk
Because of this theorem
we can write the exterior derivative of an element
of the quantum plane as
\eq
da = (\chi_c * a) V^c
\label{daflat}
\en
Thus $da$ is expressed in terms of the ``q-vielbein" $V^c$.
\sk
The value and action of $\chi_s$ on the coordinates
is easily computed, cf. the
definition in  (\ref{chiin}) :
\eq
\chi_c (x^a)=-q_{c\bu} \de^a_c~,~~~
\chi_c * x^a=-q_{c\bu} \T{a}{c}~; \label{actval}
\en
using (\ref{copex}) we find the deformed Leibniz rule
\eq
\chi_c * (ab)=(\chi_c * a) \f{c}{c} * b + a \chi_c *  
b~.\label{Leibflat}
\en
\noi From (\ref{actval}) the exterior derivative of $x^a$ is:
\eq
dx^a=-q_{c\bu} \T{a}{c} V^c=-V^c\T{a}{c} \label{dxa}
\en
\noi
[use (\ref{VTcomm})]
and gives the relation between the $q$-vielbein $V^c$
and the differentials $dx^a$.
\sk
The $x^a$ and $V^b$ $q$-commute as (cf. (\ref{Vxcomm})):
\eq
V^a x^b= q_{a\bu}~ x^b V^a
\en
\noi and via eq. (\ref{dxa}) and (\ref{PRTT33})
 we find the $dx^a, x^b$commutations :
\eq
dx^a x^b=q_{ab}~x^b dx^a
\en
\noi After acting on this equation with $d$ we obtain the
commutations between the differentials:
\eq
dx^a \we dx^b= - q_{ab}~ dx^b \we dx^a \label{wedplane}~.
\en
The commutations between the partial derivatives are given in
eq.(\ref{qLie3}).
\sk
All the relations of this section
are covariant under the
$SO_{q,r=1}(N)$ action:
\eq
x^a \longrightarrow \T{a}{b} \otimes x^b~.
\en
Notice that the partial derivatives $\chi_c$,
and in general all the tangent vectors $\chi$
and vector fields $\chi *$
of this paper have ``flat" indices. To compare $\chi_c *$
with partial derivative operators with ``curved" indices,
we need to define the operators $\stackrel{\leftarrow}{\chit_s} :$
\eq
\stackrel{\leftarrow}{\chit_s}\!(a) \equiv - {1\over q_{b\bu}}  
(\chi_b * a)
\kappa (\T{b}{s})~, \label{partialcurved}
\en
\noi so that
\eq
d\,a= \stackrel{\leftarrow}{\chit_s}\!(a) ~dx^s
\en
\noi which is equation (\ref{daflat}) in ``curved" indices.
The action of $\stackrel{\leftarrow}{\chit_s}$ on the coordinates is
\eq
\stackrel{\leftarrow}{\chit_s}\!(x^a)=\de^a_s I~,~~
\stackrel{\leftarrow}{\chit_s}\!(bx^a)=b\stackrel{\leftarrow}{\chit_s} 
\!(x^a)+
q_{sa}\stackrel{\leftarrow}{\chit_s}\!(b)\;x^a~\label{Leibpart}
\en
\sk
The tangent vector fields $\chi_c *$ of this paper and the partial  
derivatives
$\stackrel{\leftarrow}{\chit_s}$
are derivative operators that act  ``from the right to the left"
as it is seen from their deformed
Leibniz rule (\ref{Leibflat}), (\ref{Leibpart}). This explains the
inverted arrow on $\stackrel{\leftarrow}{\chit_s}$.
We can also define derivative operators acting from the left
to the right, as in ref.s \cite{qplane},  using the antipode $\kappa$  
which is
antimultiplicative.
For a generic quantum group the vectors ${-\kpm}(\chi_i) \equiv  
-\chi_i \ci
\km$
act from the left and we also have
\eq
d\,a=(\chi_i*a)\omega^i=\omega^i({-\kpm}(\chi_i)*a)
\en
as is seen  from $\kp (\chi_i)=-\chi_j \kp (\f{j}{i})$ and $\kpm  
(\f{k}{j})
\f{i}{k}
= \de^i_j$ [third line of (\ref{copf})].

We then define the partial derivatives
\eq
\chit_s (a)\equiv
{\kappa}^{-1}(\T{i}{s})\;\kpm(\chi_i)*a~,\label{partialcurvedr}
\en
so that
\eq
d\,a=dx^s\;\chit_s(a)~.
\en
The action  of $\chit_s$ on the coordinates is
\eq
\chit_s (x^a)=\de^a_s I~,~~
\chit_s({x^a b})=\chit_s(x^a)~b+q_{as}x^a\chit_s(b)
\en
\sk
{}From eqs. (\ref{partialcurved}), (\ref{partialcurvedr})
and (\ref{qLie3}), or directly from  $d^2=0$ and $dx^a \we dx^b=dx^a  
\otimes
dx^b - q_{ab}~ dx^b \otimes dx^a$
[a consequence of (\ref{wedplane})],
one finds the following commutations between the ``curved" partial  
derivatives:
\eq
\stackrel{\leftarrow}{\chit_r}\;  \stackrel{\leftarrow}{\chit_s} =  
q_{sr}
\stackrel{\leftarrow}{\chit_s}\;
\stackrel{\leftarrow}{\chit_r}~~,~~~
{\chit_r}\;  {\chit_s} = q_{rs} {\chit_s}\;  {\chit_r}~.
\en

Finally, we note that the transformation
\eqa
& &\xi^a={1\over \sqrt{2}}(x^a+x^{a'}),~~~\sma{$a \leq n$}\\
& &\xi^{n+1}={i\over \sqrt{2}}(x^n - x^{n+1})\\
& &\xi^a={1\over \sqrt{2}}(-x^a+x^{a'}),~~~\sma{$a > n+1$}
\ena
\noi defines real coordinates $\xi^a$ for the even dimensional
orthogonal quantum plane
$Fun_{q,r=1}(ISO(n+1,n-1)/SO(n+1,n-1))$
endowed with the conjugation ii)
discussed in Section 2. Moreover on this basis the metric becomes
diagonal.
Likewise it is possible to define antihermitian $\chi$ and real
$\Omega, V$.
\sk
For $n=2$ the results of this section immediately yield the
bicovariant calculus on the $q$-Minkowski
space, i.e. the
multiparametric orthogonal quantum plane
$Fun_{q,r=1}(ISO(3,1)/SO(3,1))$.
 \vskip 1cm


\app{ The Hopf algebra axioms}


A Hopf algebra over the field $K$ is a
unital algebra over $K$ endowed
with the  linear maps:
\eq
\D~: ~~A \rightarrow A\otimes A,~~
\epsi~:~~ A\rightarrow K,~~
\kappa~:~~A\rightarrow A
\en
\noi satisfying the following
properties $\forall a,b \in A$:
\eq
(\D \otimes id)\D(a) = (id \otimes \D)\D(a) \label{axiom1}
\en
\eq
(\epsi \otimes id)\D(a)=(id \otimes \epsi)\D(a) =a \label{axiom2}
\en
\eq
m(\kappa\otimes id)\D(a)=m(id \otimes \kappa)\D(a)
=\epsi(a)I \label{kappadelta}
\en
\eq
\D(ab)=\D(a)\D(b)~;~~\D(I)=I\otimes I
\en
\eq
\epsi(ab)=\epsi(a)\epsi(b)~;~~\epsi(I)=1
\en
\noi where $m$ is the multiplication map $m(a\otimes b)
= ab$.
{} From  these axioms we deduce:
\eq
\kappa(ab)=\kappa(b)\kappa(a)~;~~\D[\kappa(a)]=\tau(\kappa\otimes
\kappa)\D(a)~;~~\epsi[\kappa(a)]=\epsi(a)~;~~\kappa(I)=I
\en
where $\tau(a \otimes b)=b\otimes a$ is the twist map.

\vskip 1cm


\vfill\eject
\end{document}